\documentclass[a4paper,12pt]{article}
\usepackage[british]{babel}
\usepackage{amsfonts}
\usepackage{amsmath}
\title{Direct Construction of Aperiodic Tilings with the Hat Monotile}
\author{Ulrich Reitebuch \\ \small{ulrich.reitebuch@fu-berlin.de}}
\usepackage{graphicx}
\begin{document}
\newcommand{\ww}{0.33\textwidth}
\newcommand{\www}{0.28\textwidth}
\newcommand{\wwww}{0.22\textwidth}
\maketitle
\begin{abstract}
In 2023, the quest for an aperiodic monotile was answered by the ``hat'' monotile. In this article, structures in this aperiodic tiling are discovered, which allow for a direct computation of the tiling, similar to well-known methods for the Penrose tilings.
\end{abstract}

\section{Introduction}
The "hat" monotile was introduced in 2023 by David Smith, Joseph Samuel Myers, Craig S. Kaplan, and Chaim Goodman-Strauss \cite{2023monotile}.
Intersecting a regular triangle mesh $U$ and its dual hexagon mesh, the Euclidean plane is split into congruent kites. The hat monotile consists of eight of these kites, and there are two possible hats, one hat and its mirror image.
\begin{figure}[bth]
	\begin{center}
		\includegraphics[width=\www]{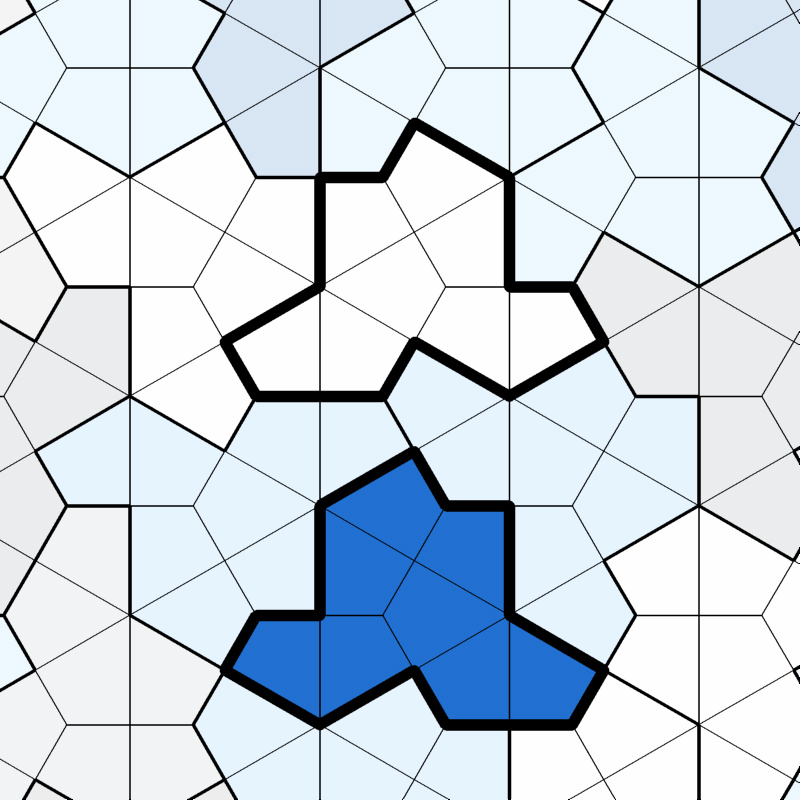}
		\hspace{0.05\textwidth}
		\includegraphics[width=\www]{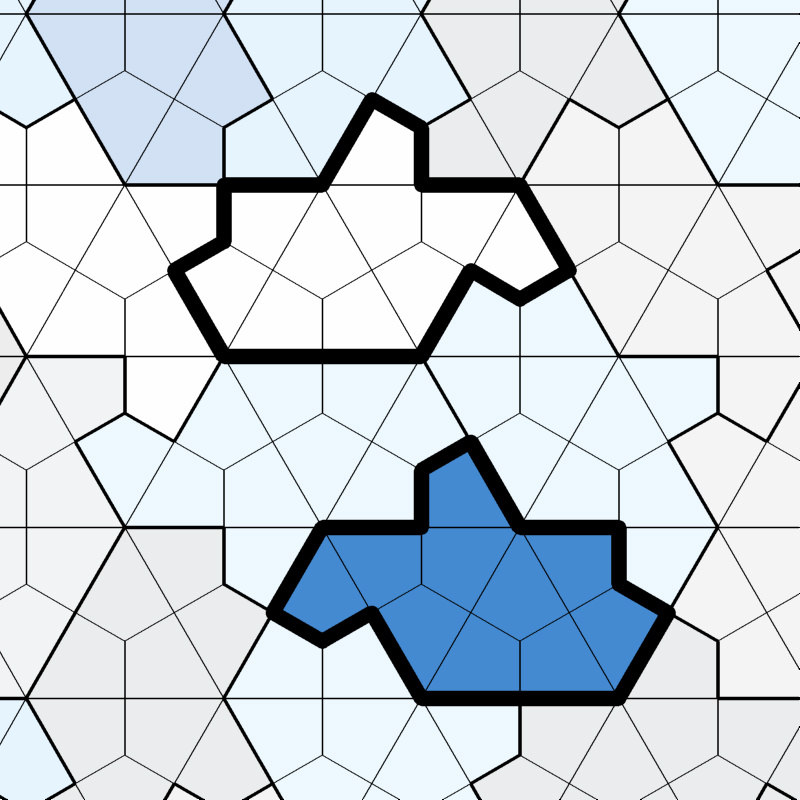}
		\caption{Flipped and undlipped monotile in the tilinigs with $8$--kites hat monotiles and $10$--kites monotiles.}
		\label{flippedAndUnflipped}
	\end{center}
\end{figure}

Both types of hats are necessary for tiling the plane, but not in the same frequency. We will refer to the more frequent hat type in a tiling as ``unflipped'', to the less frequent as ``flipped''.
The roles of flipped and unflipped hat type can be exchanged.
The hat monotile is in \cite{2023monotile} presented as a one--parameter family of monotiles, we will use the classical hat monotile, described by parameters $(1, \sqrt{3})$ and the ten-kite tile described by parameters $(\sqrt{3}, 1)$, as shown in Figure~\ref{flippedAndUnflipped}.

We will use for the golden ratio
\begin{center}
$\Phi = \frac{\sqrt{5}+1}{2}$, \\
$\varphi = \frac{\sqrt{5}-1}{2}$, 
\end{center}
with the usual golden ratio formulas
\begin{center}
$\Phi^{2} = \Phi+1$, \\
$\varphi^{2} = 1 - \varphi$, \\
$\Phi \varphi = 1$.
\end{center}

\begin{figure}[bt]
	\begin{center}
		\includegraphics[width=\www]{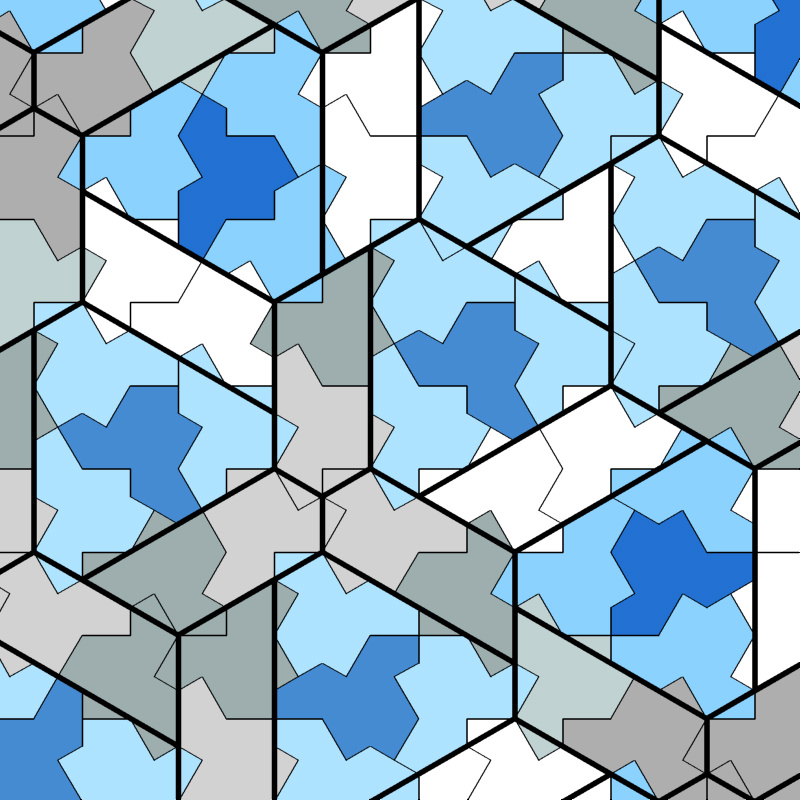}
		\caption{metatiles highlighting structure in the hat tiling.}
		\label{metaTiles}
	\end{center}
\end{figure}

Hat tiles are drawn in \cite{2023monotile} as dark blue hats (these are the flipped hats), unflipped hats are light blue, white, and grey hats.
These colours are related to the tiling by metatiles in section $2$ - a ``H'' metatile represents a flipped and three light blue unflipped hats, a ``T'' metatile represents an isolated white hat, a ``P'' metatile represents a pair of white hats, and a ``F'' metatile represents a pair of grey hats (F metatiles are the third part of a ``fylfot'' structure), see Figure~\ref{metaTiles}.

\section{Dual Triangulation}
\begin{figure}[tb]
	\begin{center}
		\includegraphics[width=\wwww]{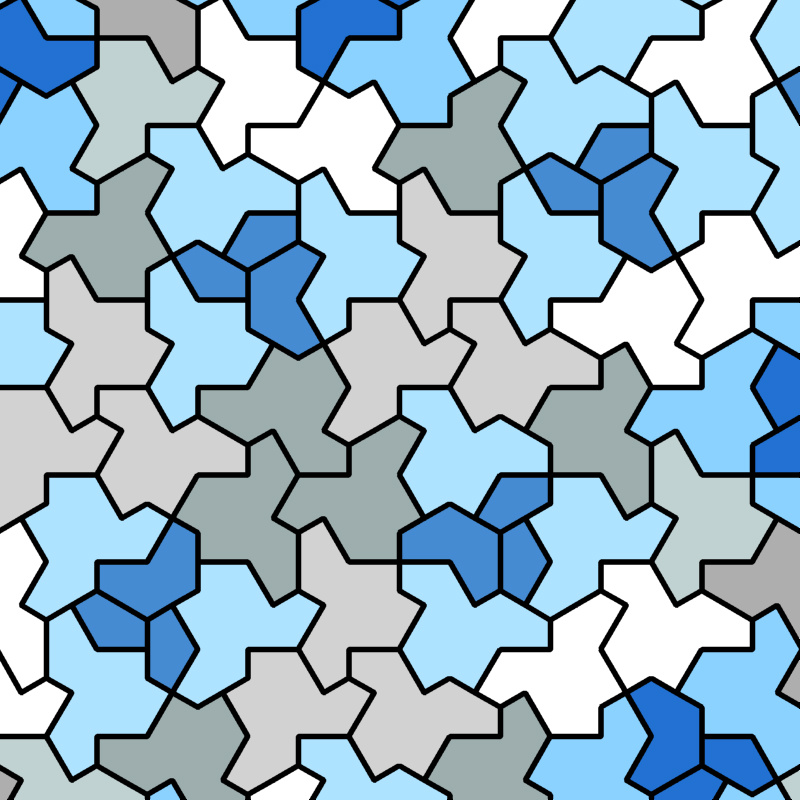}
		\hspace{0.01\textwidth}
		\includegraphics[width=\wwww]{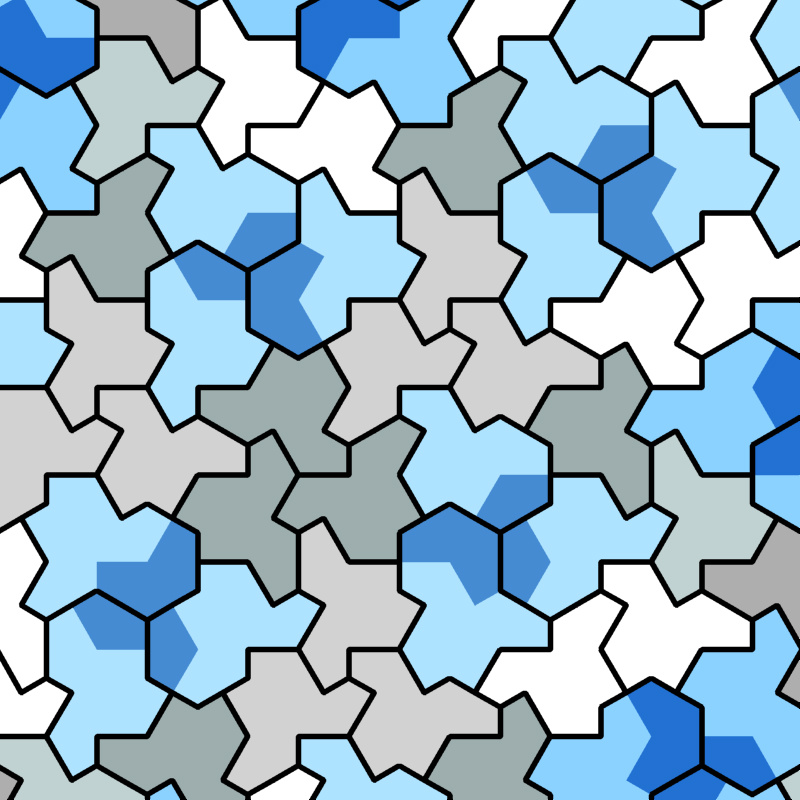}
		\hspace{0.01\textwidth}
		\includegraphics[width=\wwww]{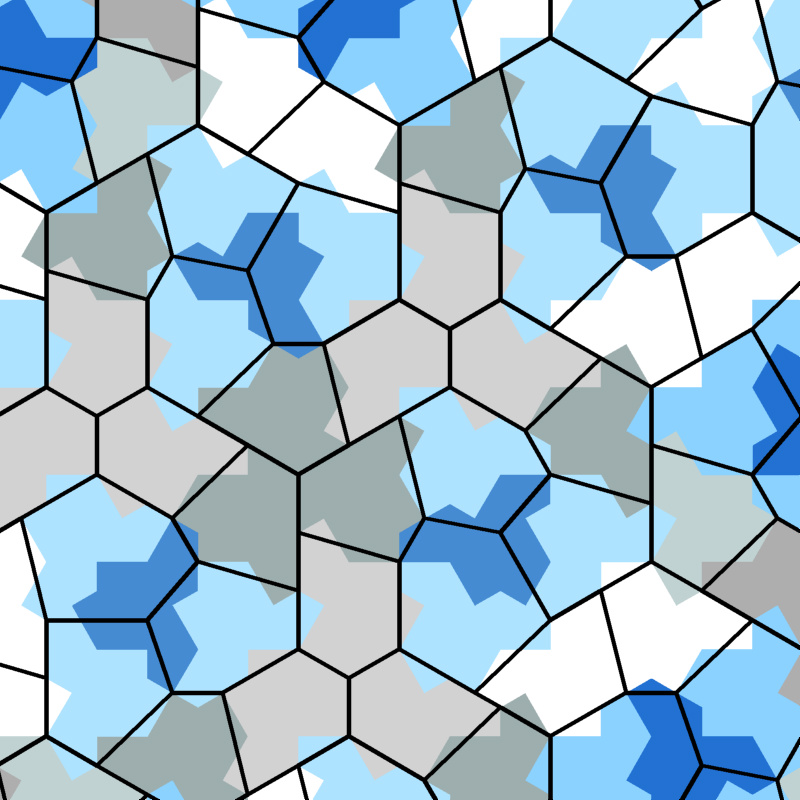}
		\hspace{0.01\textwidth}
		\includegraphics[width=\wwww]{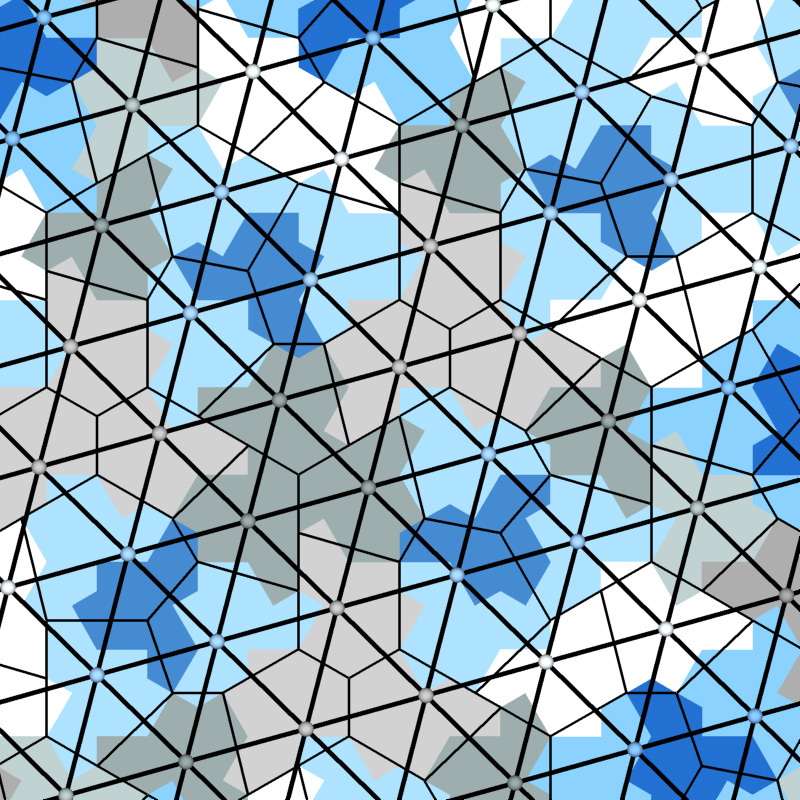}
		\caption{Left to right: Split flipped hats into three parts -- Join parts to adjacent light blue hats -- Hexagon mesh connectivity -- Dual triangle mesh.}
		\label{dualTriangles}
	\end{center}
\end{figure}

It was shown that the tiling obtained by joining each flipped hat with a certain (unflipped) neighbour is related to a regular hexagon grid (compare Figure 2.2 in \cite{2023monotile}).
If each of the flipped hats is split into three parts along kite edges, and these three parts are joined to the three adjacent light blue hats, the new tiles will even already have the connectivity of a regular hexagon grid, see Figure~\ref{dualTriangles}.

In this article, We will use its dual grid, a regular triangle grid $T$, to construct any aperiodic hat tiling directly, using only two parameters (similar to the pentagrid method by N. G. De Bruijn \cite{1981algebraicTheory} \cite{1996RemarksOnPenrose}).
At first, $T$ is combinatorial information, but we will consider several useful geometric realisations.
Each vertex of the triangle grid will represent one unflipped hat -- when these all are constructed, the remaining gaps will have the shape of flipped hats.

\section{Triangle Grid Lines}
\begin{figure}[tb]
	\begin{center}
		\includegraphics[width=\www]{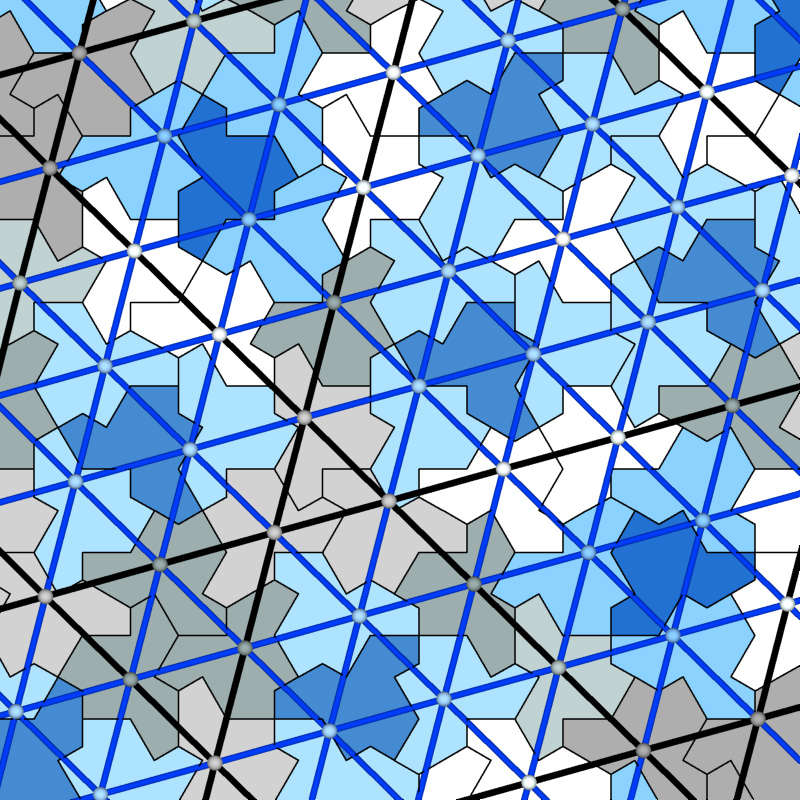}
		\hspace{0.05\textwidth}	
		\includegraphics[width=\www]{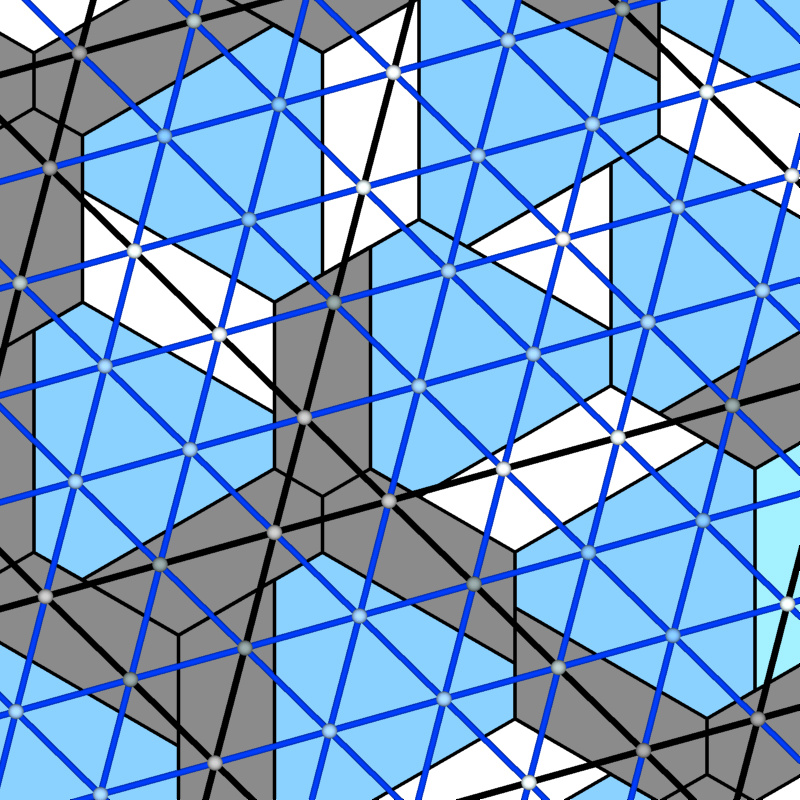}
		\hspace{0.05\textwidth}	
		\includegraphics[width=\www]{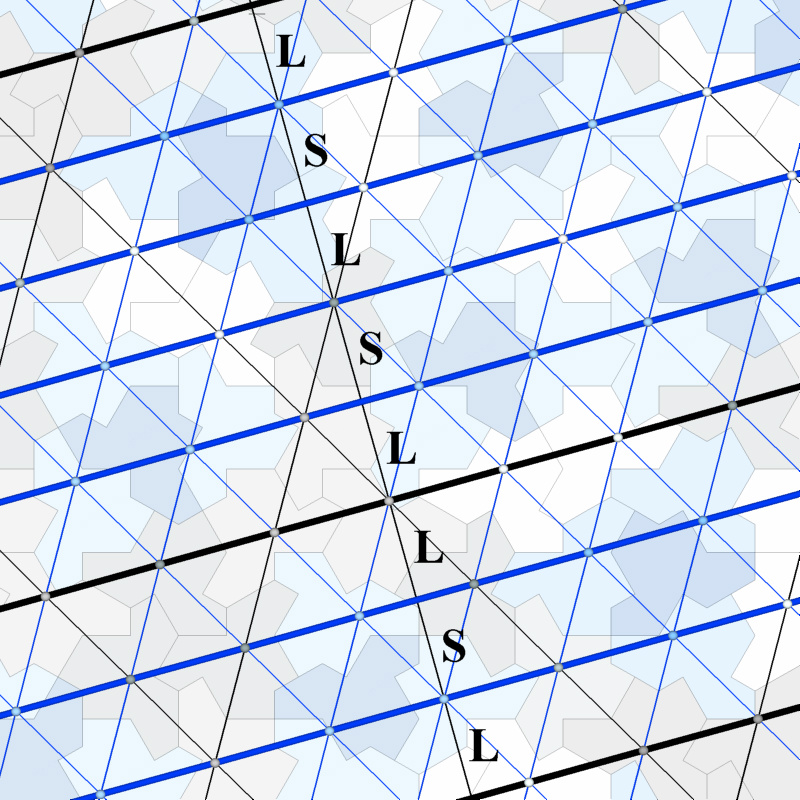}
		\caption{Left to right: Black and blue triangulation lines in the hat tiling -- Black and blue triangulation lines in the metatiles -- Labeling gaps between lines.}
		\label{blackAndBlue}
	\end{center}
\end{figure}

Looking at the grid lines in the triangulation $T$, we can see that there are lines which pass only through vertices representing grey and white hats (we will refer to these as ``black'' lines), and there are other lines which pass through vertices representing all kinds of hats (we will refer to these as ``blue'' lines), see Figure~\ref{blackAndBlue}.
The black lines roughly follow the direction of P and F metatiles.

Looking at parallel grid lines, the black lines are always next to two blue lines, and blue lines always come in pairs.
Between two black lines, there may be one or two pairs of blue lines.
Note, there are neighbouring blue lines that are not considered a ``pair'', in each blue pair at least one of the two blue lines should have one black neighbour line.

There are different types of vertices in the triangle grid: vertices crossed by three blue lines represent either light blue hats or isolated white hats belonging to a ``T'' metatile, vertices crossed by one black and two blue lines represent either a white hat belonging to a ``P'' metatile or a grey hat, and vertices crossed by two black lines and one blue line represent grey hats.
There are no vertices in the triangulation $T$ which are crossed by three black lines.

The colouring of one set of parallel lines looks very much like the patterns in the Amman bars in the Penrose tilings:
taking the distance between a pair of blue lines as the short distance $S$ and all distances between neighbouring lines which do not belong to the same pair of blue lines (this might be two blue or a black and a blue line) as long distance $L$, this looks like the line pattern we get by the substitution rule for Ammann bars in Penrose tilings:
\begin{center}
$S \rightarrow L$, \\
$L \rightarrow S+L$.
\end{center}
We will refer to this pattern of $L$ and $S$ distances as ``Fibonacci pattern''.

\subsection{Golden Ratio}
To show that the black and blue lines in $T$ in fact do follow the Fibonacci pattern, we will use the converged shape of the metatiles: 
there is a substitution rule for metatiles given in \cite{2023monotile}, but by the substitution, the tiles change their shape - angles are preserved, but the ratio between different edge length is changed. There is however also a converged shape given for the metatiles, where the substitution preserves the shape - but the converged metatiles do not exactly fit to the hat tiling.
The ratios between edge lengths of these converged tiles can be expressed by terms containing the golden ratio, compare Figure~2.8 in \cite{2023monotile}.

\begin{figure}[bt]
	\begin{center}
		\includegraphics[width=\wwww]{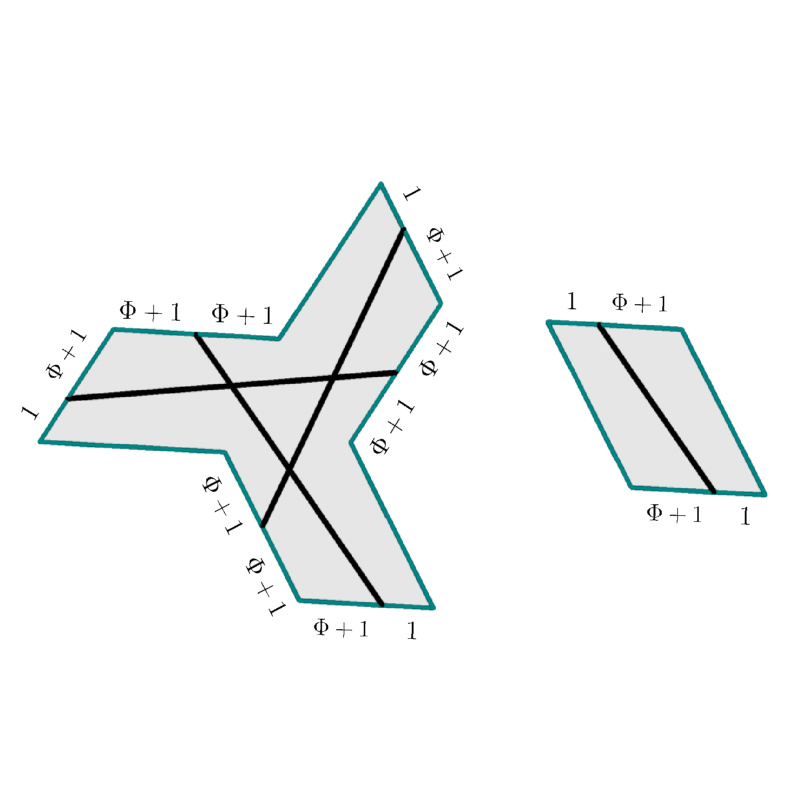}
		\hspace{0.01\textwidth}
		\includegraphics[width=\wwww]{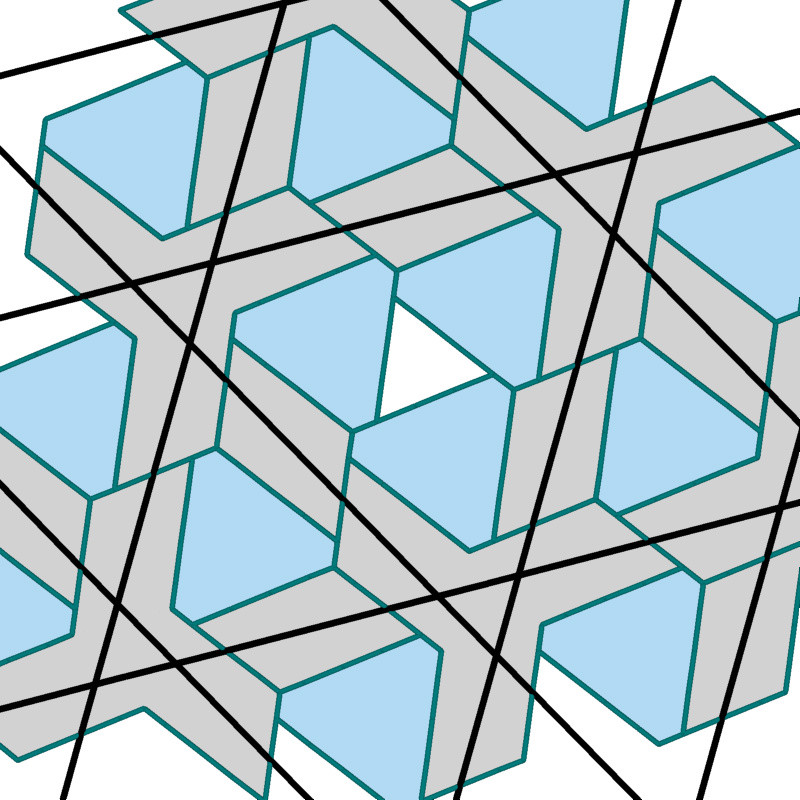}
		\hspace{0.01\textwidth}
		\includegraphics[width=\wwww]{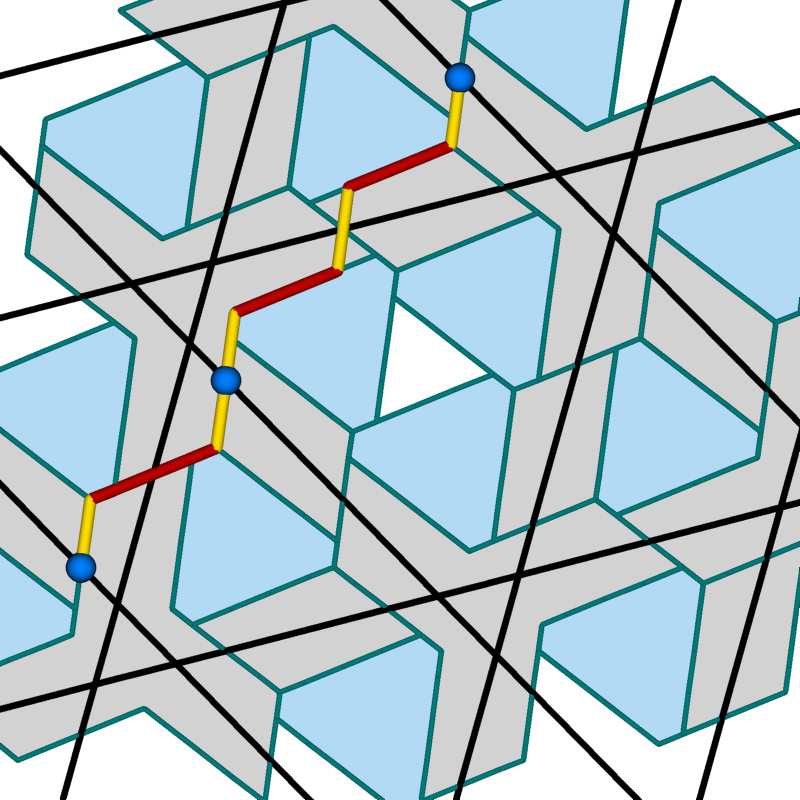}
		\hspace{0.01\textwidth}
		\includegraphics[width=\wwww]{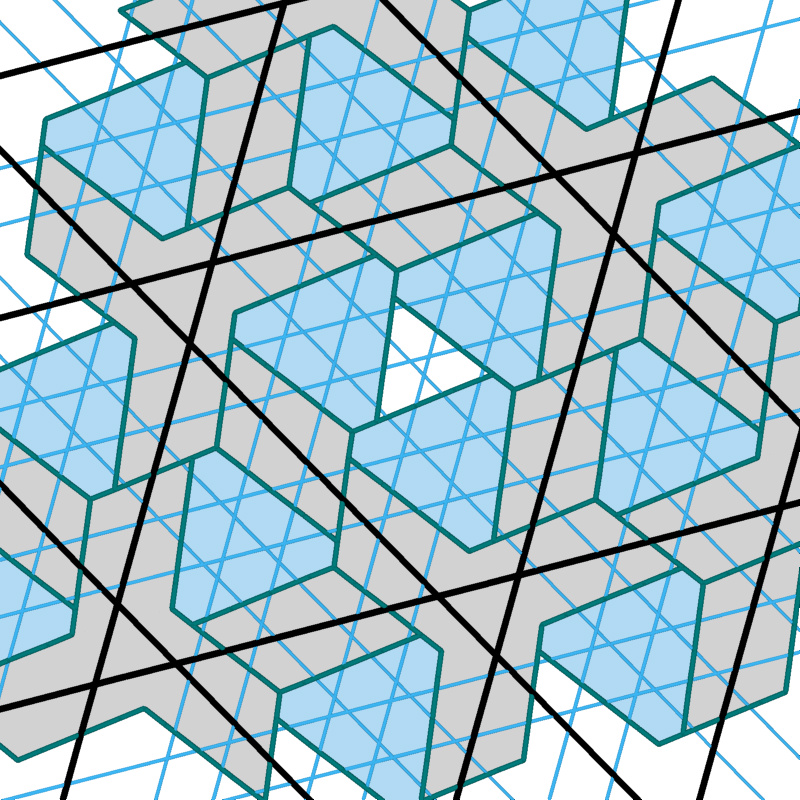}
		\caption{Left to right: Black lines on F and P converged metatiles -- Straight black lines in the converged metatiles -- Red and yellow parts sum up to $2+2\Phi$ in the short distance, to $2+4\Phi = \Phi(2+2\Phi)$ in the long distance between blue points -- blue lines painted between black lines in Fibonacci pattern.}		
		\label{convergedMetatiles}
	\end{center}
\end{figure}

The black lines of the triangulation $T$ can be drawn on the converged metatiles as straight lines, similar to the Ammann bars in the Penrose tilings.
A suitable position for black lines is shown in Figure~\ref{convergedMetatiles}. By this realisation, black lines completely stay inside ``P'' and ``F'' metatiles.
Painting only these black lines on the converged monotiles, there are long and short distances to the next parallel black line, and the ratio between short and long distances between black lines is $\varphi$.
This can be directly computed from the edge lengths of the metatiles, see Figure~\ref{convergedMetatiles}.

The blue lines can be drawn in between the black lines, such that distances between a blue pair of lines is $\varphi$ times the distance between neighbouring lines that are not a blue pair.
however, painting black and blue lines with these distances, the lines in three directions never intersect in one point -- there are only intersections of two lines.
Thus the vertices of the triangulation $T$ are not drawn as one intersection point of three lines, but as three close-by intersection points of only two lines each.
So, for observations on the substitution rule and the golden ratio in the line pattern we will use only the black lines and omit the blue lines - blue lines could be inserted at every substitution level, replacing long $L$ distances between black lines by $LSLSL$ and short $S$ distances between black lines by $LSL$, and then in the fine pattern colouring all lines next to a $S$ symbol blue, lines between two $L$ symbols stay black.
\begin{figure}[tb]
	\begin{center}
		\includegraphics[width=\www]{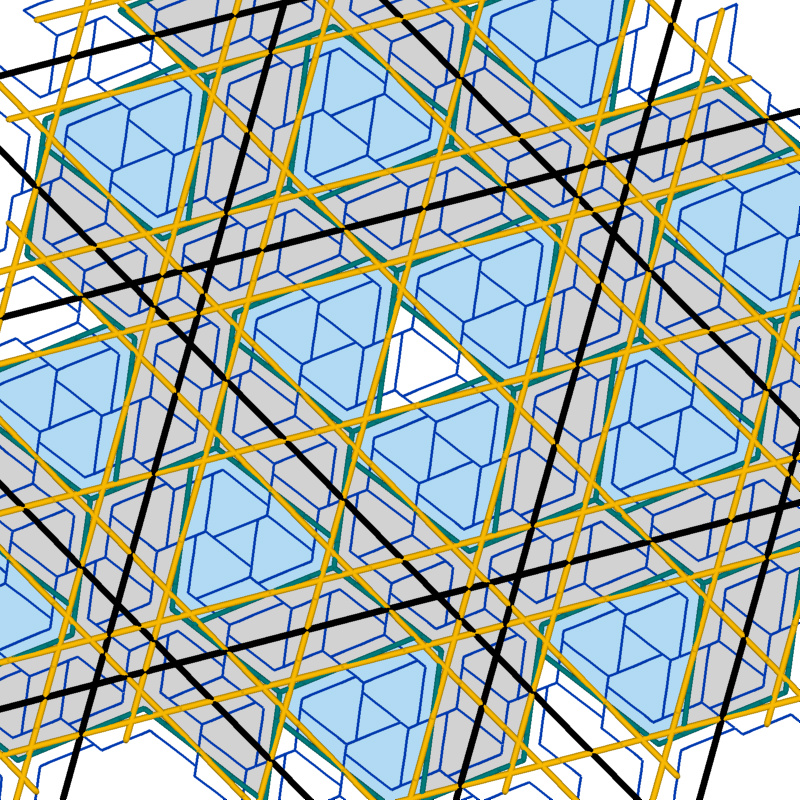}
		\hspace{0.01\textwidth}
		\includegraphics[width=\www]{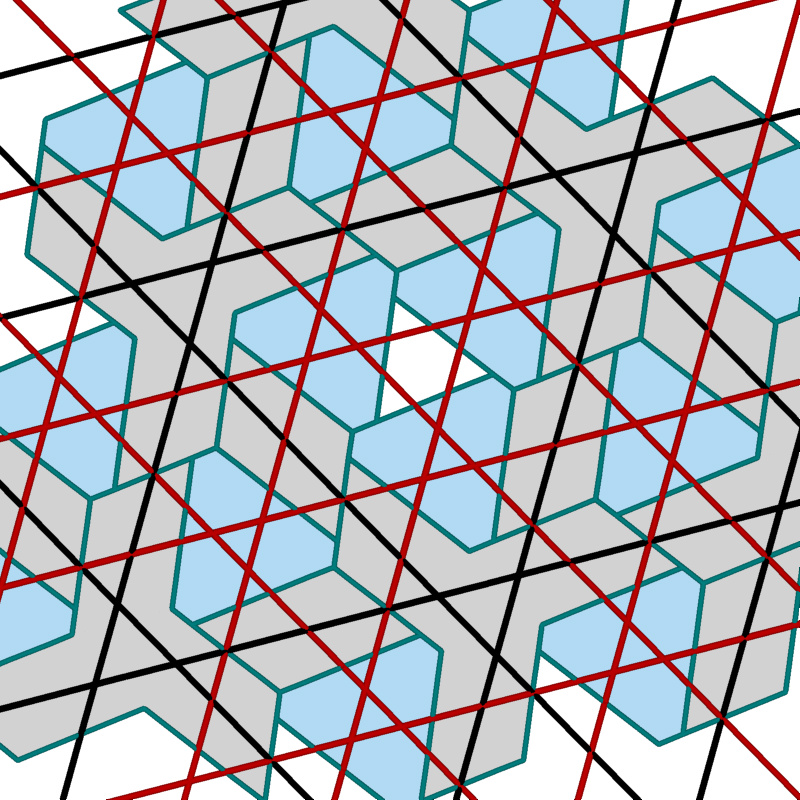}
		\hspace{0.01\textwidth}
		\includegraphics[width=\www]{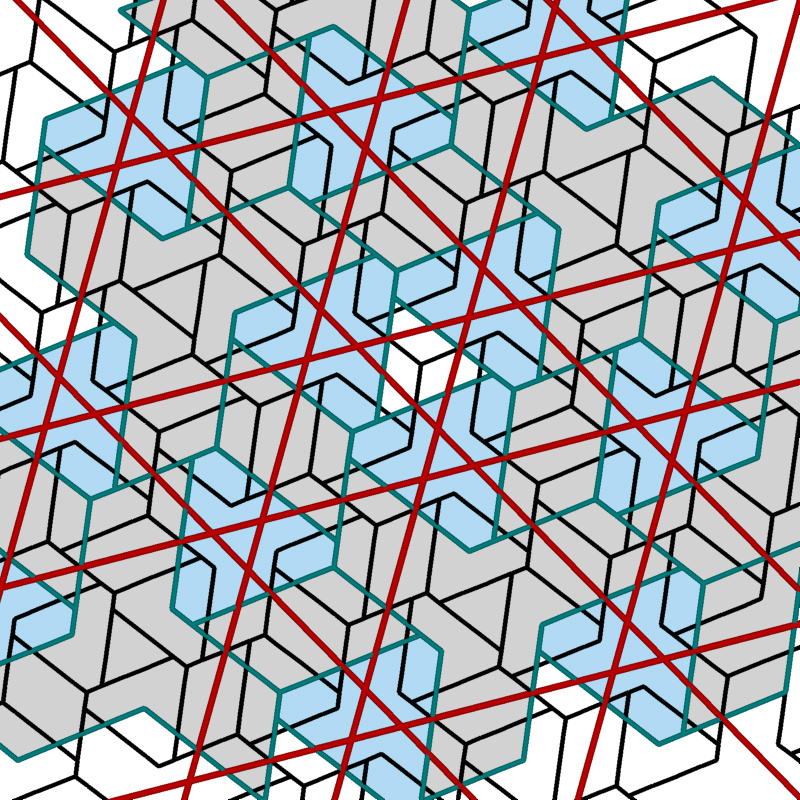}
		\caption{Left to right: Metatile substitution rule and finer lines as in \cite{2023monotile} -- Slower substitution for lines -- Resulting slower metatile substitution.}
		\label{metatileSubstitution}
	\end{center}
\end{figure}

The substitution for metatiles shown in Figure~2.8 in \cite{2023monotile} actually performs two steps of the Fibonacci line substitution, so by the substitution we get for the $S$ and $L$ distances between black lines 
\begin{center}
$S \rightarrow S+L$, \\$L \rightarrow S+2L$.
\end{center}
The geometric realisation of the black lines shown in Figure~\ref{convergedMetatiles} is not only creating more lines with smaller $L$ and $S$ distances, but the ratio of different distances in the substituted lines is again exactly $\varphi$.

For only one step of line substitution, we get a slower substitution for the metatiles, see Figure~\ref{metatileSubstitution}.
By the substitution rule for the metatiles, the short and long distances between black lines follow the Fibonacci substitution rule, as the Ammann bars in the Penrose tiling do.
Also for the slower substitution, the new lines after performing the substitution will again have small and large distances with ratio $\varphi$.

These observations on the Fibonacci pattern, which holds individually for each of the three line directions, can also be used to prove aperiodicity of the tiling: for any translational symmetry of the tiling, the ratio between $S$ and $L$ labels for at least two directions of lines would have to be rational, but $\varphi$ is an irrational number.


\section{6D Index Vectors}
\begin{figure}[bt]
	\begin{center}
		\includegraphics[width=\wwww]{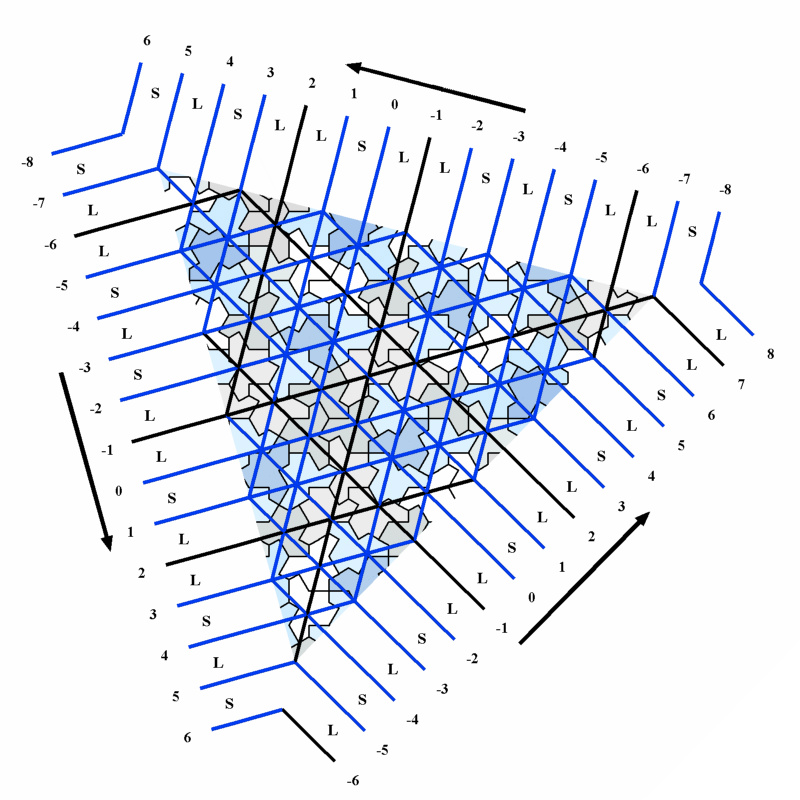}
		\hspace{0.01\textwidth}
		\includegraphics[width=\wwww]{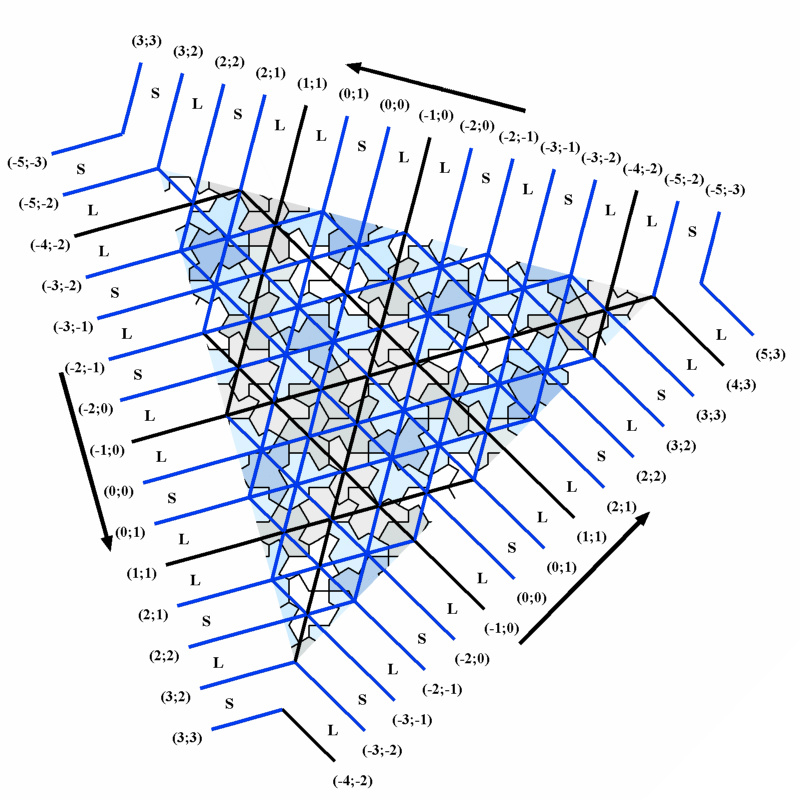}
		\hspace{0.01\textwidth}
		\includegraphics[width=\wwww]{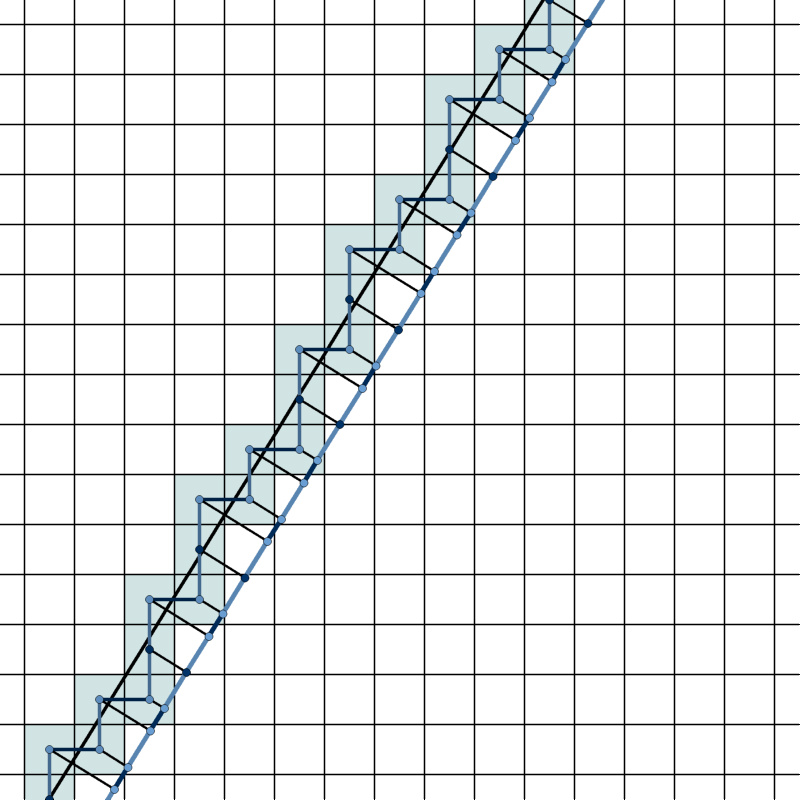}
		\hspace{0.01\textwidth}
		\includegraphics[width=\wwww]{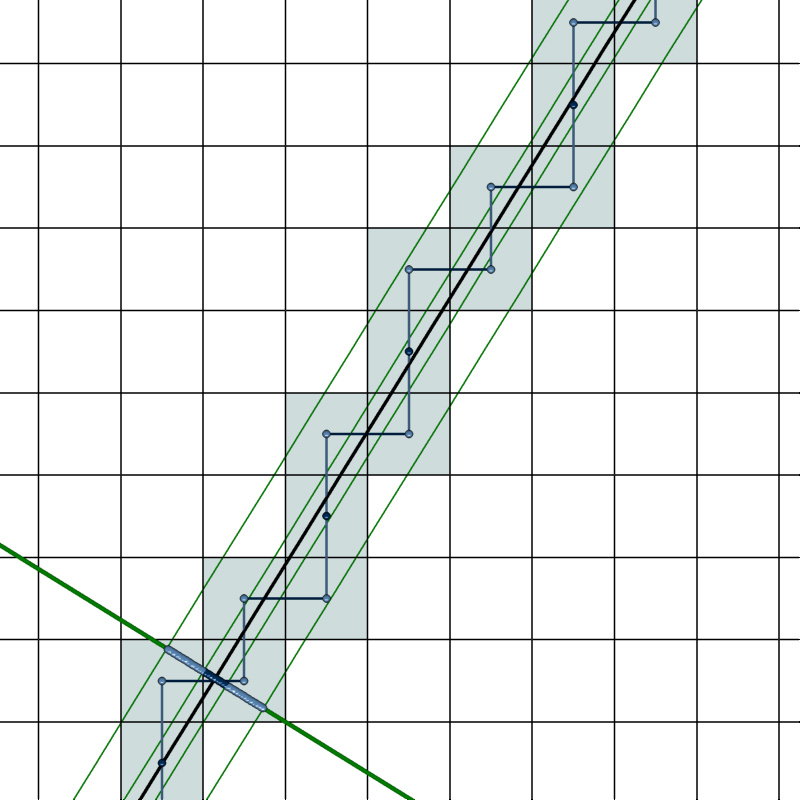}
		\caption{Left to right: Integer labels for triangulation lines -- Fibonacci indices for lines -- Projection of square centres to a line -- Orthogonal projection.
		}	
		\label{lineLabels}
	\end{center}
\end{figure}

Now we create an integer index vector for each vertex of the triangulation $T$.
First, we label the parallel lines in each of the three directions $k~\in~\{0, 1, 2\}$ by integer indices, in ascending order in three directions with angle $\frac{2\pi}{3}$ in between, as shown in Figure~\ref{lineLabels}. So each vertex $v_i$ of $T$ is assigned three integers $n_0(i),n_1(i),n_2(i)$ according to the indices of the three lines passing through $v_i$.

In a second step, we assign $S$ and $L$ distances in order of the Fibonacci pattern to each of the three sets of parallel lines in the $T$.
This can be done by the projection method for each of the three directions $k~\in~\{0, 1, 2\}$ individually: for three given real numbers $d_{k}$, we define for each index $n_{k}(i)$ two indices
	\begin{align}
		a_{k}(i) &= \lfloor \varphi n_{k}(i) + d_{k}\rfloor, \label{a}\\
		b_{k}(i) &= n_{k}(i)-a_{k}(i). \label{b}
	\end{align}
This gives a $6D$ integer index vector $(a_0(i), b_0(i), a_1(i), b_1(i), a_2(i), b_2(i))$ for each vertex $v_i$ of the triangulation.

A well-known geometric interpretation is shown in Figure~\ref{lineLabels}:
when a straight line $g_k(x)=\Phi x + d_k$ is drawn into a unit square grid with vertices in $\mathbb{Z}^2$, projecting the centres of all squares intersected by the line $g_k$ to a parallel line creates a Fibonacci pattern of short ($S$) and long ($L$) distances.
By the indices $a_{k}$ and $b_{k}$ for two neighbouring projected points $v_i$ and $v_j$, if $a_{k}(i) = a_{k}(j)$ the distance is labelled by $S$, otherwise it is labelled by $L$.

Projecting the same square centre points to a line orthogonal to $g_k$, all points are projected to a finite segment of the orthogonal line, see Figure~\ref{lineLabels}.
This orthogonal projection also gives information on the types of points: 
points between two $L$ symbols are projected to the centre part of the segment, points next to a $S$ and a $L$ Symbol are projected closer to one of the ends of the segment, dividing the segment in three sections with ratios $\varphi^2 : \varphi^3 : \varphi^2$.
This is similar to the pentagonal and pentagram-shaped regions used to classify vertices of the Penrose rhombus tiling in De Bruijn's work \cite{1981algebraicTheory}.
We will see that there is a similar kind of classification for the unflipped hat tiles.

Lines in the triangulation are coloured black and blue as above: if a line is next to one distance labelled as $S$, it is coloured blue, if both adjacent distances are labelled as $L$, the line is coloured black.
To avoid three black lines through one vertex, we choose the values $d_{k}$ such that $d_0 + d_1 + d_2 = 0$.
Furthermore, the $d_{k}$ should be chosen such that the line $g_k$ does not pass through any grid point for $k \in \{0, 1, 2\}$ to avoid ambiguities, i.e. $d_{k}$ cannot be written as $i+\varphi j$ with $i, j \in \mathbb{Z}$.

By the choice of the $d_k$ parameters, a Fibonacci pattern for each line direction is created, and reversely, for each Fibonacci pattern there exists a suitable value of $d_k$, so these parameters are one--to--one with the possible hat tilings.

\section{3D Coordinates from 6D Integer Indices}
There are several possibilities to compute $3D$ or $2D$ positions from the $6D$ index vectors for the vertices of the triangulation $T$. 
\begin{figure}[bt]
	\begin{center}
		\includegraphics[width=\wwww]{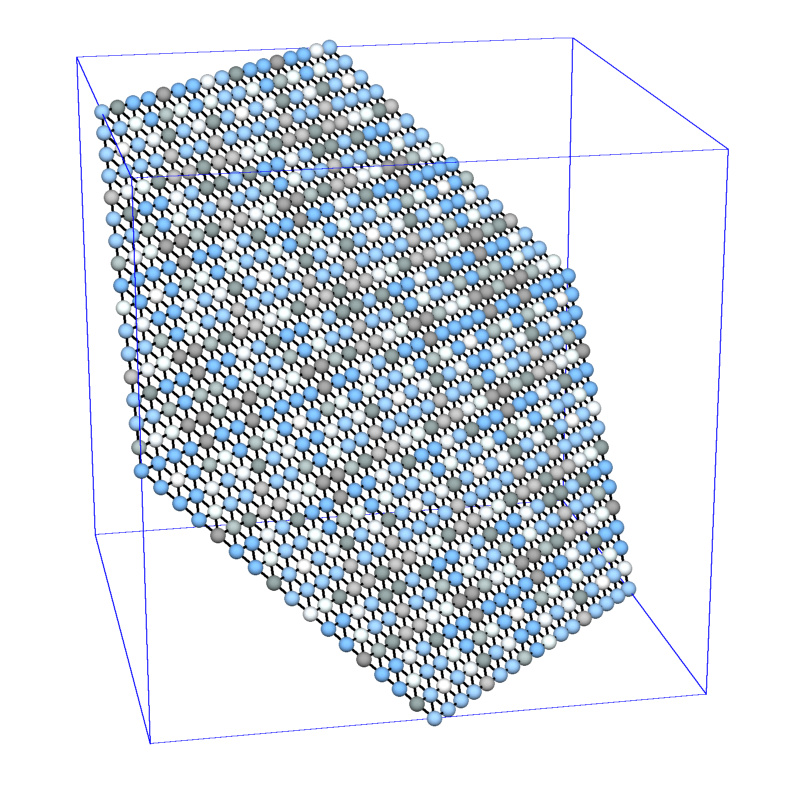}
		\hspace{0.01\textwidth}
		\includegraphics[width=\wwww]{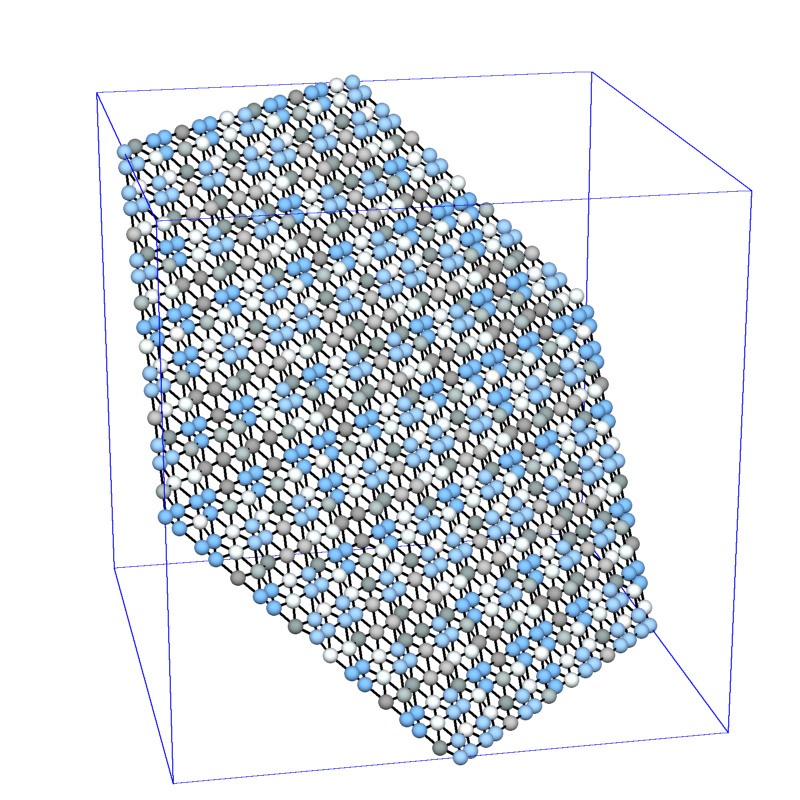}
		\hspace{0.01\textwidth}
		\includegraphics[width=\wwww]{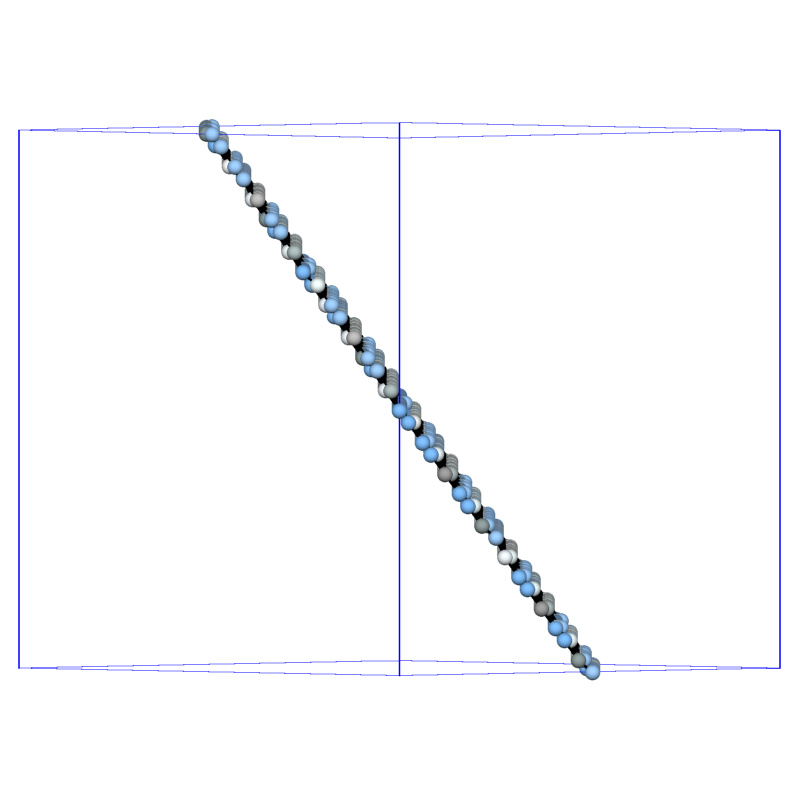}
		\hspace{0.01\textwidth}
		\includegraphics[width=\wwww]{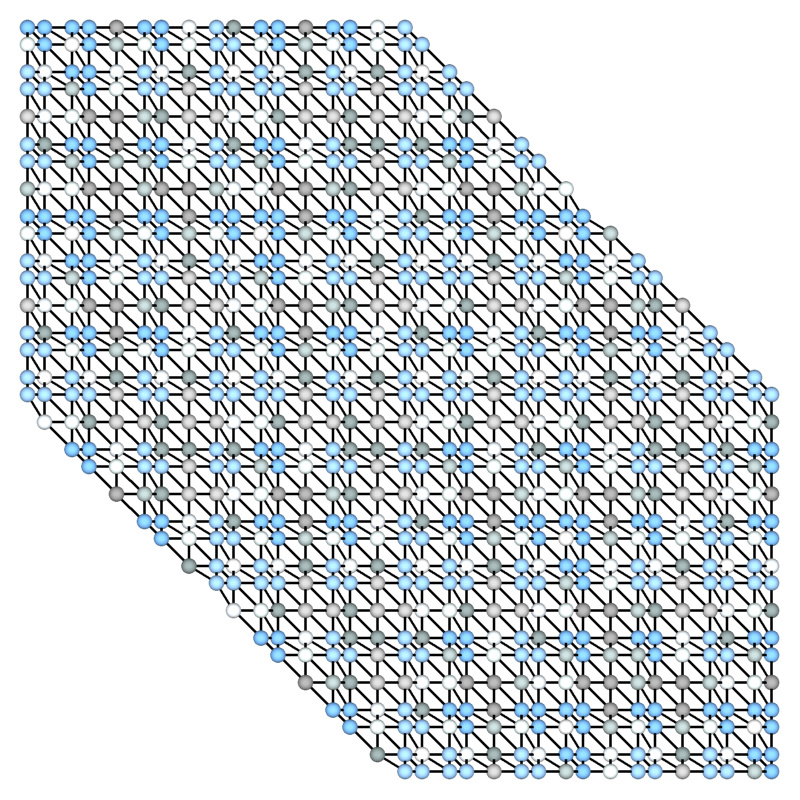}
		\caption{Left to right: $3D$ cooridnates $v_r$ -- $3D$ cooridnates $v_f$ -- $v_f$ allignment on two planes - orthogonal projection of $v_f$ to the x-y-plane.
		}	
		\label{plane3D}
	\end{center}
\end{figure}
By 
\begin{center}
$
v_r(i) = \begin{pmatrix} x_r(i) \\ y_r(i) \\ z_r(i) \end{pmatrix} = \begin{pmatrix}a_0(i)+b_0(i) \\ a_1(i)+b_1(i) \\ a_2(i)+b_2(i)\end{pmatrix}
$
\end{center}
the vertices of the triangulation $T$ are mapped to a plane perpendicular to the vector
\begin{center}
$n=\begin{pmatrix} 1 \\ 1 \\ 1 \end{pmatrix}$,
\end{center}
and in this case the triangulation consists of geometrically regular triangles.

If we instead use 
\begin{center}
$
v_f(i) = \begin{pmatrix} x_f(i) \\ y_f(i) \\ z_f(i) \end{pmatrix} = \begin{pmatrix}a_0(i) \Phi + b_0(i) \\ a_1(i) \Phi + b_1(i) \\ a_2(i) \Phi + b_2(i)\end{pmatrix}
$,
\end{center}
we increase each $3D$ coordinate by $1$ for short distances and by $\Phi$ for long distances, and the vertices of the triangulation align on two parallel planes, both perpendicular to the same vector $n$, see Figure~\ref{plane3D}; vertices are coloured by the types of unflipped hats they represent.

\begin{figure}[bt]
	\begin{center}
		\includegraphics[width=\www]{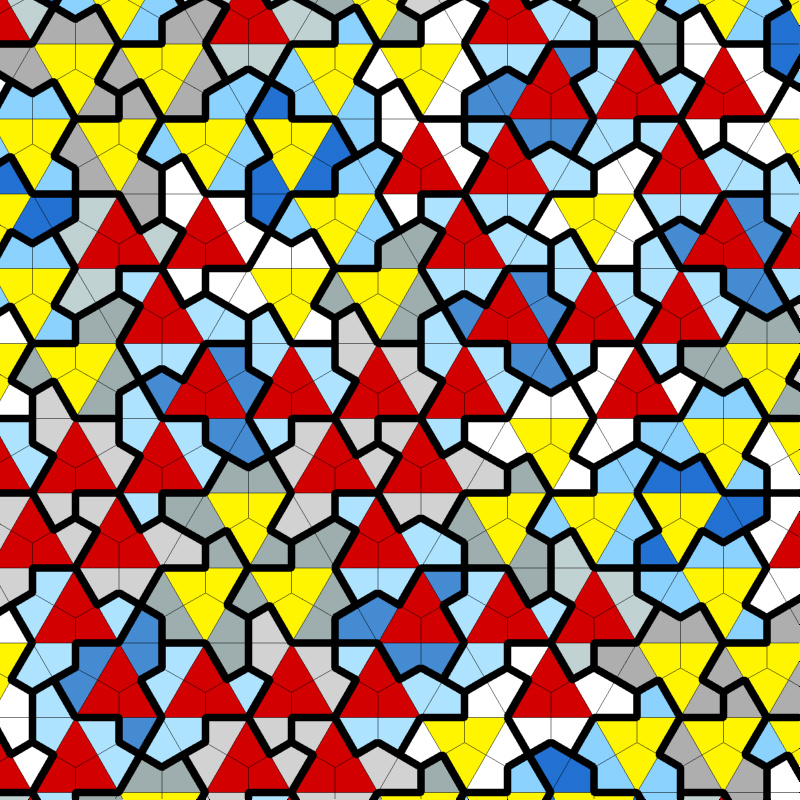}
		\hspace{0.01\textwidth}
		\includegraphics[width=\www]{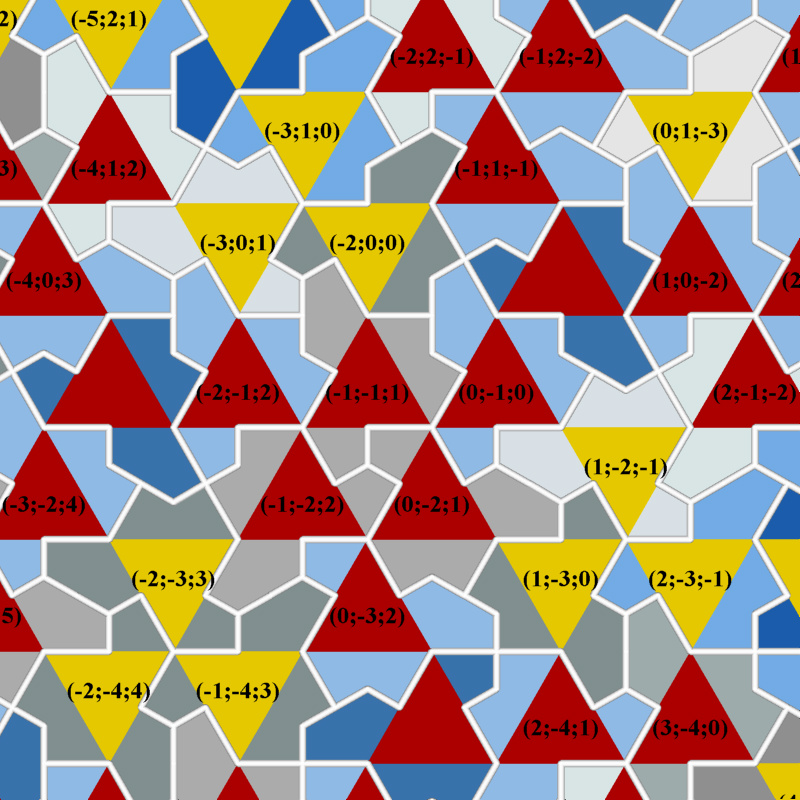}
		\hspace{0.01\textwidth}
		\includegraphics[width=\www]{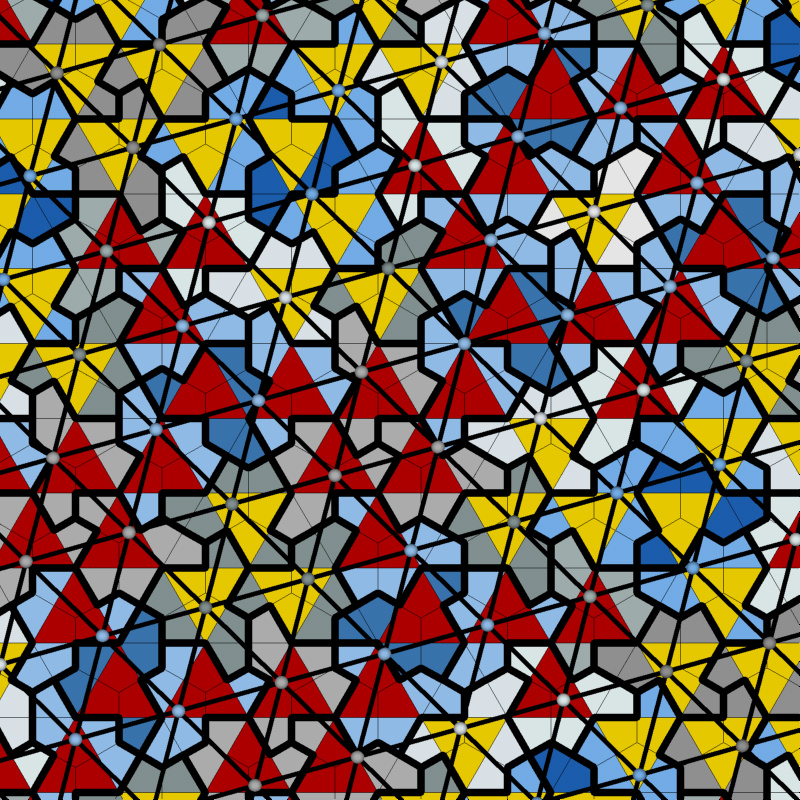}
		\caption{Left to right: Centre triangles in all hats, triangles pointing up in red, triangles pointing down in yellow -- $c_k$ labels for centre triangles in unflipped hats -- dual triangulation $T$ with vertices in centre triangles.
		}	
		\label{centerTriangles}
	\end{center}
\end{figure}

Each hat contains one complete regular triangle from the initial triangulation $U$ (we will call this the ``centre triangle'' of the hat), built of three of its eight kites, see Figure~\ref{centerTriangles}.
There are two classes of centre triangles, some are pointing up, some are pointing down.
neighbouring light blue centre triangles always are in the same class (and in the same class as the centre triangle of the flipped hat in the middle of three of the light blue hats).
Neighbouring white hats always belong to different classes.
In each fylfot structure, the centre triangles of the three grey hats at the centre are in the same class, the other three are in the other class.

From the $6D$ indices for the vertices of the triangulation, we can compute a 3D index vector 
\begin{align}
c(i) = \begin{pmatrix} c_0(i) \\ c_1(i) \\ c_2(i) \end{pmatrix} = \begin{pmatrix}a_0(i)+b_0(i)-b_1(i) \\ a_1(i)+b_1(i)-b_2(i) \\ a_2(i)+b_2(i)-b_0(i)\end{pmatrix},\label{c}
\end{align}
which describes the position of the centre triangle of hat $i$ -- each of the three indices counts gaps in one direction of parallel lines in the unit triangulation from which the kites were generated.

The triangulation $T$ representing the unflipped hats can also be embedded in $2D$, such that its vertices are contained in the centre triangles of the hats they represent.
To achieve this, the triangulation $T$ has to be drawn with edge length $2 \Phi$ and angle $\alpha$ relative to $U$, where $\tan(\alpha)=\frac{\sqrt{3}}{3+2\Phi}$.

Note, for the reflected hats, exchanging the roles of (more frequent) flipped and (less frequent) unflipped hats, the index vector $c(i)$ has to be computed as
\begin{align}
c(i) = \begin{pmatrix} c_0(i) \\ c_1(i) \\ c_2(i) \end{pmatrix} = \begin{pmatrix}a_0(i)+b_0(i)-b_2(i) \\ a_1(i)+b_1(i)-b_0(i) \\ a_2(i)+b_2(i)-b_1(i)\end{pmatrix},\label{c2}
\end{align}
and for the embedded dual triangulation $T$ the sign of the rotation angle $\alpha$ has to be flipped.

\begin{figure}[bt]
	\begin{center}
		\includegraphics[width=\wwww]{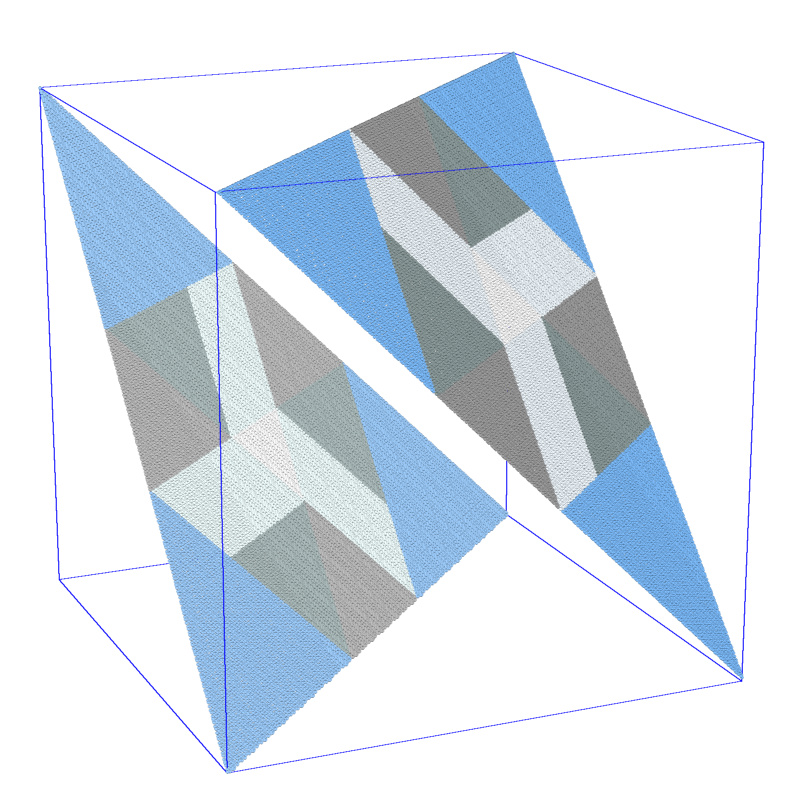}
		\hspace{0.01\textwidth}
		\includegraphics[width=\wwww]{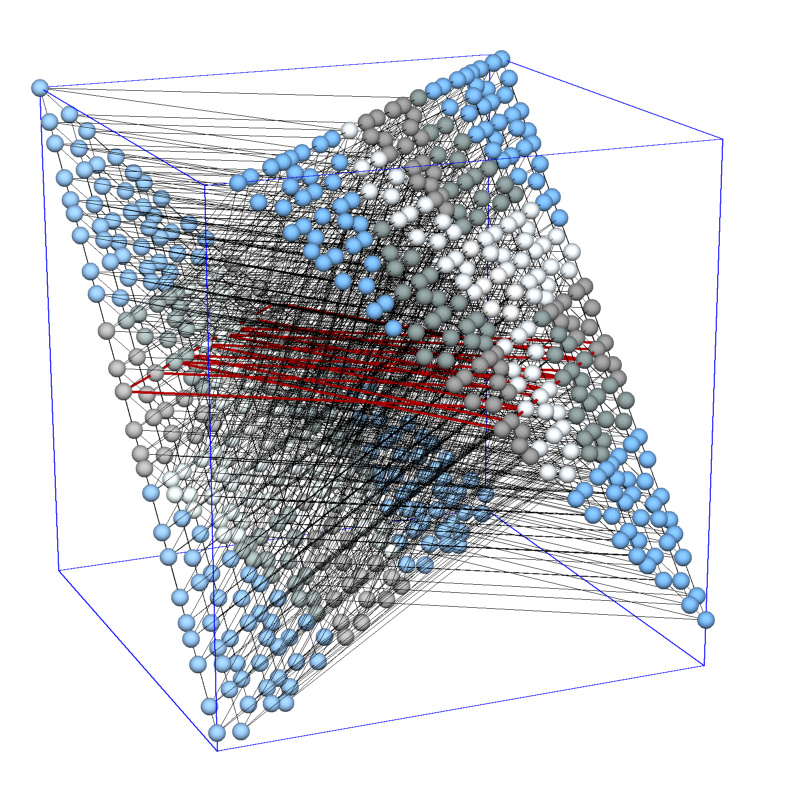}
		\hspace{0.01\textwidth}
		\includegraphics[width=\wwww]{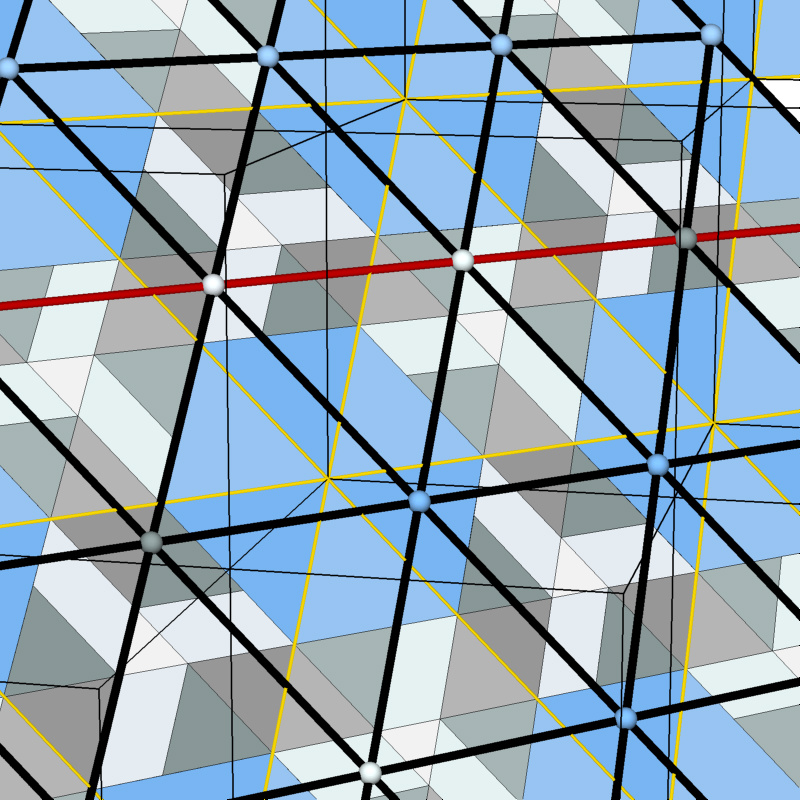}
		\hspace{0.01\textwidth}
		\includegraphics[width=\wwww]{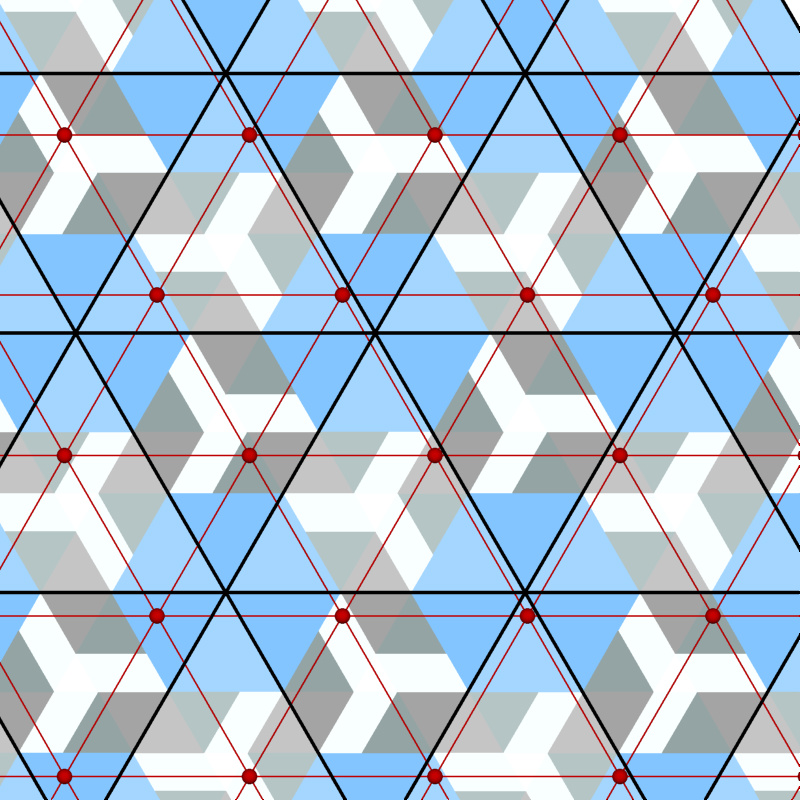}
		\caption{Left to right: $v_b(i)$ coordinates, vertices coloured by hat types - one triangulation line highlighted in red, folding back in the box -- repeating the box, the triangulation can be unfolded to a regular triangulation in the colouring pattern -- hat colour pattern in 2D, grid size is $\Phi$ times the grid size of $T$ (embedded as regular triangle grid in 2D, shown in red).
		}	
		\label{box}
	\end{center}
\end{figure}

As in the projection of square centres to a line orthogonal to the generating line $g_k$ shown in Figure~\ref{lineLabels} in, we can compute a different geometric realisation $T_b$ of $T$, which will indicate information about the hats represented by the vertices of $T$:
with 
\begin{center}
$
v_b(i) = \begin{pmatrix} x(i) \\ y(i) \\ z(i) \end{pmatrix} = \begin{pmatrix}a_0(i) - b_0(i) \Phi \\ a_1(i) - b_1(i) \Phi \\ a_2(i) - b_2(i) \Phi \end{pmatrix}
$,
\end{center}
we increase the coordinate by $1$ for every long distance but decrease the coordinate by $\Phi$ for every short distance; in this case all vertices of $T$ are mapped into a cube with edge length $1+\Phi$.
The vertices are all positioned on two triangles spanned by six cube corners, for both these triangles $n$ is a normal vector, see Figure~\ref{box}.
One of the triangles contains all vertices representing hats with centre triangle pointing up, the other contains all vertices representing hats which centre triangles pointing down.

If we now colour the vertices of $T_b$ according to the colours of the hats they represent, we can observe an interesting pattern:
vertices close to the three corners of the cube (if all coordinates in $3D$ are closer than $1$ to the respective coordinate of the closest corner of the cube) represent light blue hats, there are six rhombi containing all vertices representing grey hats, and the remaining area contains all vertices representing white hats.
The two centre triangles in the white area contains the vertices representing isolated white hats, the six white parallelograms contain vertices representing pairs of white hats.

When we now again take copies of this $1+\Phi$ cube and stack those, the copies of the two triangles can be joined to an infinite periodic pattern $P$ of regular $\sqrt{2} (1+\Phi)$ triangles, and the embedding of the triangulation becomes a regular triangulation $T_p$ on the triangle pattern $P$ - the vertices of the triangulation now have coordinates $v_{\Phi} = \Phi v_r$, and the relative coordinates in the triangle pattern $P$ can be used to get the types of the represented hats.
The triangle edges in $T_p$ are $\varphi$ times as long as the triangle edges in $P$.

Instead of the colour pattern in $P$ to get hat types, the mirror image of the pattern would work as well - in that case, the role of flipped and unflipped hats has to be changed, i.e. the more frequent hat type and the less frequent hat type are both replaces by their mirror images.

So, to construct a hat tiling, we start by computing $6D$ integer coordinates for the vertices representing unflipped hats, as defined by equations~(\ref{a}) and (\ref{b}). These are for the triangulation uniquely defined, once two of the parameters $d_0$, $d_1$ and $d_2$ are chosen.
These integer coordinates can be used to compute the geometric realisations of the dual triangulation: 
From the $c_k(i)$ coordinates computed by equation~(\ref{c}) or (\ref{c2}) we know the position of the centre triangle, and from $v_{\Phi}(i)$ or $v_b(i)$ coordinates we can see the type of represented hat.
Centre triangles for flipped dark blue hats can now be already inserted at unused triangles of $U$ that share vertices with three different light blue hat centre triangles.
\begin{figure}[bt]
	\begin{center}
		\includegraphics[width=\wwww]{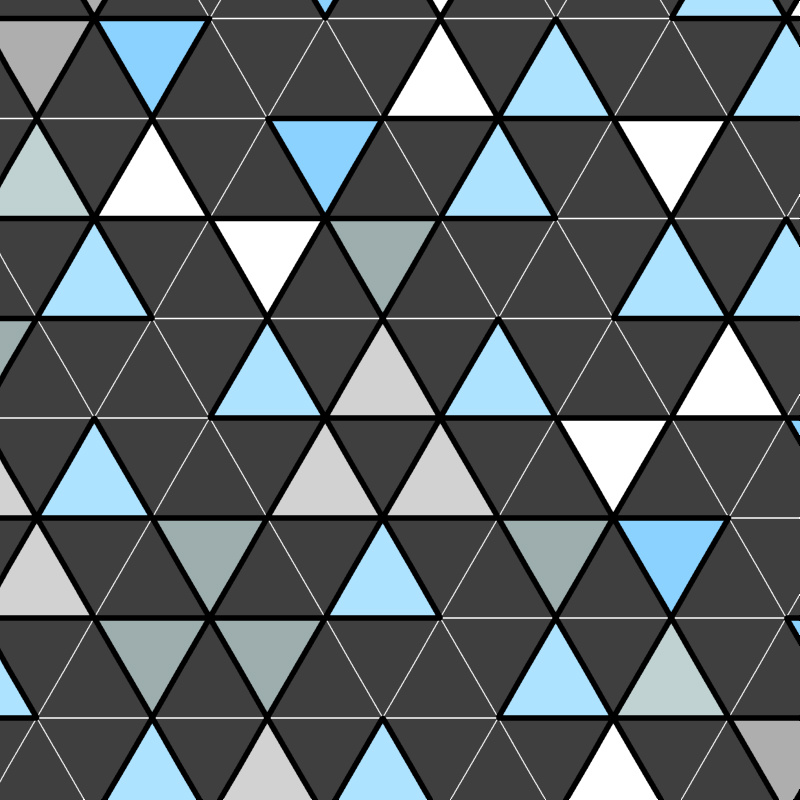}
		\hspace{0.01\textwidth}
		\includegraphics[width=\wwww]{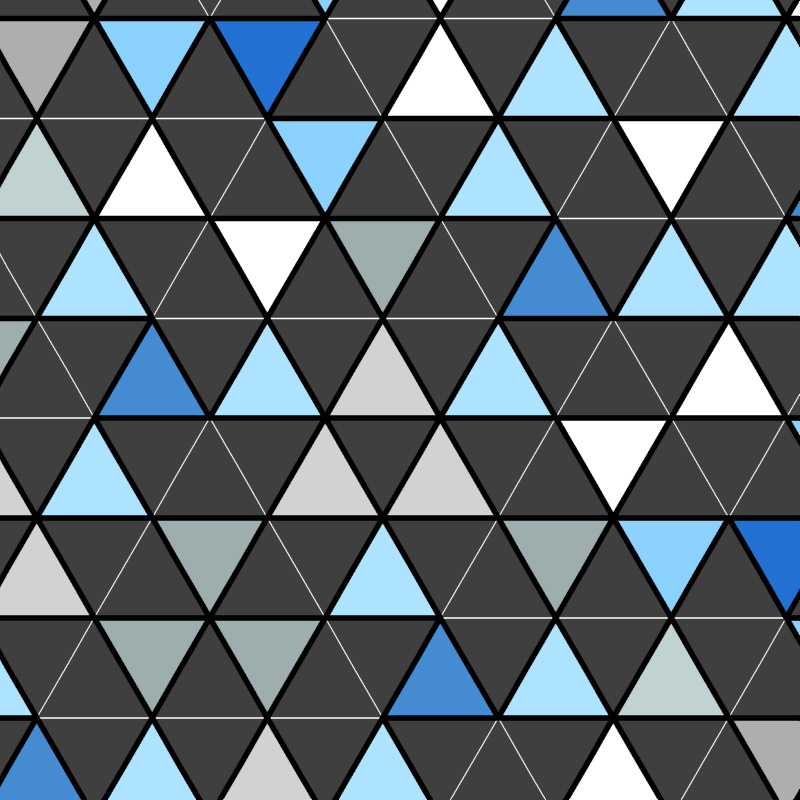}
		\hspace{0.01\textwidth}
		\includegraphics[width=\wwww]{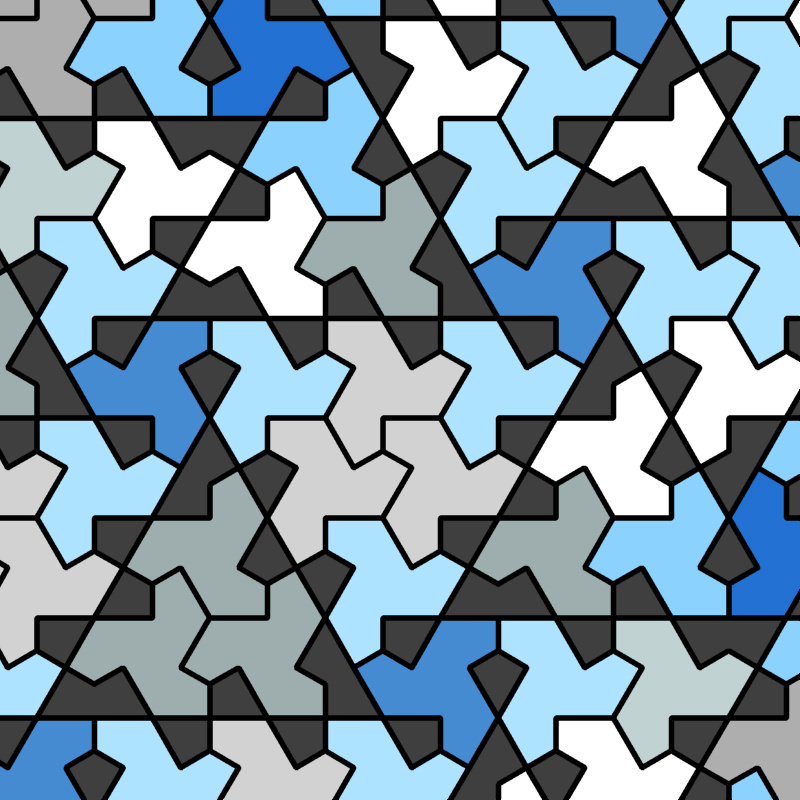}
		\hspace{0.01\textwidth}
		\includegraphics[width=\wwww]{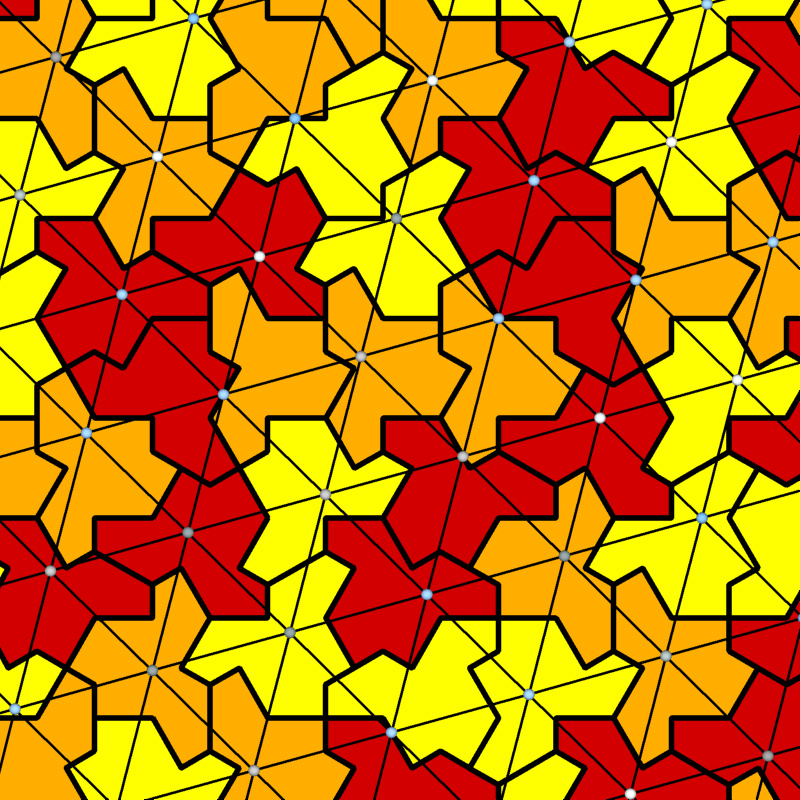}
		\caption{Left to right: Unflipped centre triangles, coloured by hat types -- Fill in centre triangles of flipped hats -- six kites for each hat, rotationally symmetric by $\frac{2\pi}{3}$ -- All complete hats coloured by orientation.
		}	
		\label{hatParts}
	\end{center}
\end{figure}

By this information we can already draw a six-kites part of each hat, consisting of the centre triangle and three cyclically attached kites, see Figure~\ref{hatParts}; only these six kites per hat would leave three possible orientations for each hat, since the the two missing kites for each hat could be attached at every corner of the centre triangle -- in some cases, not all three hats are possible, because some of the kites are already used by another (partial) hat tile.

\section{Computing Hat Orientations}
In the following observations, we will colour hats sometimes by their orientation: There are six orientations for flipped and unflipped hats -- we use three colours, such that hats with orientation rotated by $\pi$ will have the same colour.

Once the first six kites of each hat are drawn, we can see immediately, that the three grey hats at the centre of a ``fylfot'' structure have only one possible choice, since the other two choices would need kites already belonging to a different hat.
As soon as these grey hats are completed, the remaining grey hats also only have one possible orientation each.
\begin{figure}[bt]
	\begin{center}
		\includegraphics[width=\wwww]{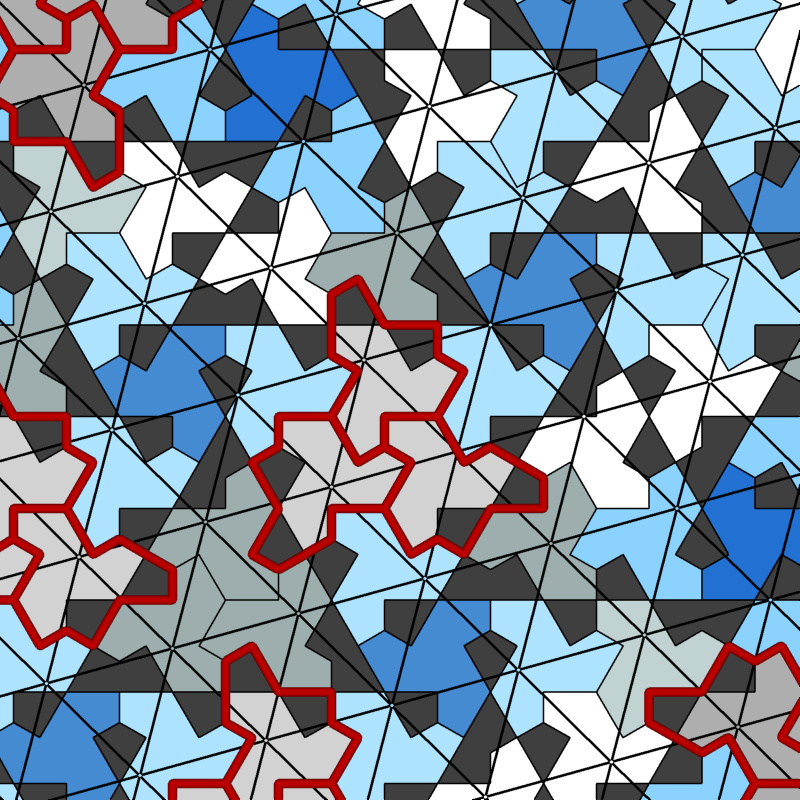}
		\hspace{0.01\textwidth}
		\includegraphics[width=\wwww]{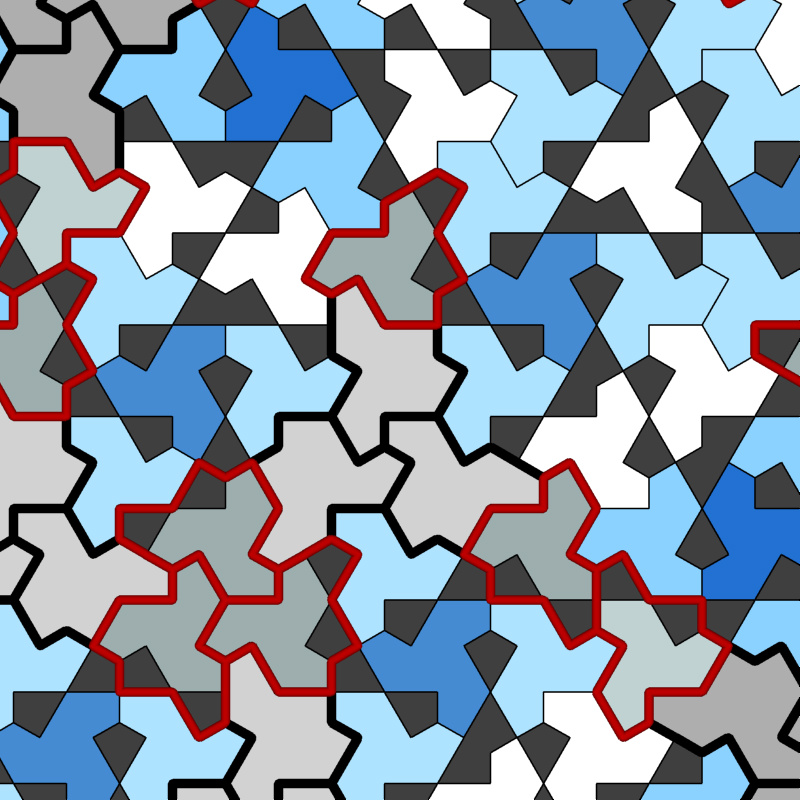}
		\hspace{0.01\textwidth}
		\includegraphics[width=\wwww]{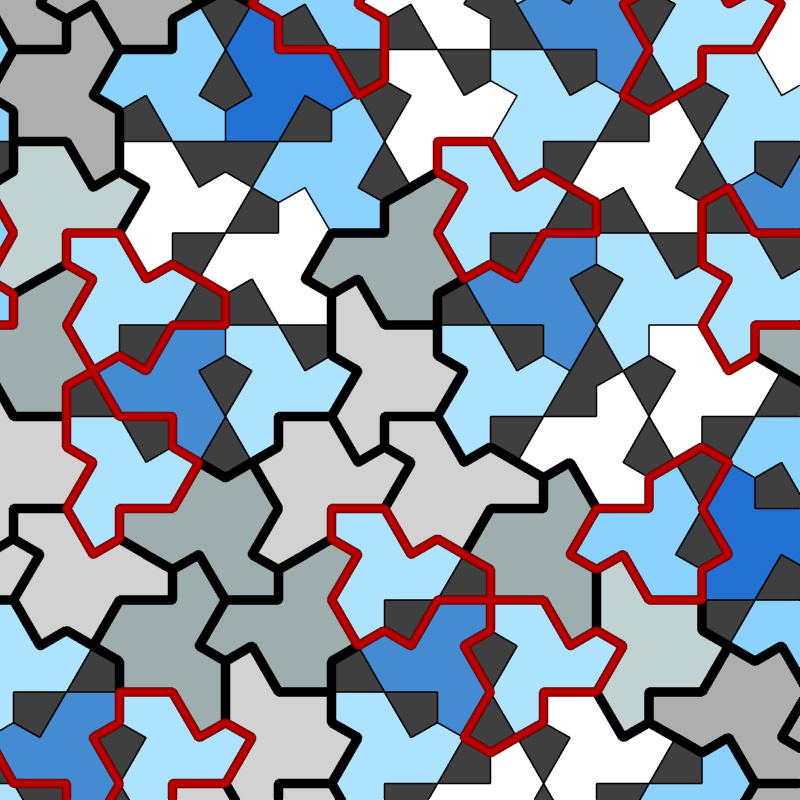}
		\hspace{0.01\textwidth}
		\includegraphics[width=\wwww]{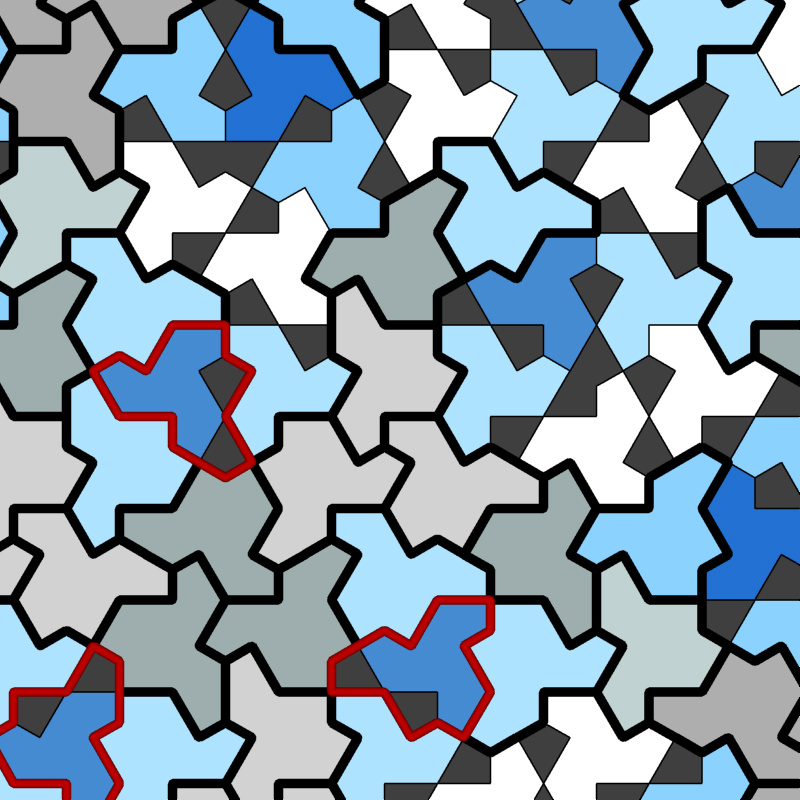}
		\caption{Orientations can be found successively. Left to right: Fylfot centre grey hats have only one possible orientation - Fill in first grey hats, remaining grey hats have only one possible orientation - fill in second grey hats, first light blue hats have only one possible orientation - fill first light blue hats, further hats have only one possible orientation.
		}	
		\label{startFindOrientation}
	\end{center}
\end{figure}
Each of the grey hats completed in the second step again leaves for one of its neighboured light blue hats only one possible orientation, see Figure~\ref{startFindOrientation}.

In a similar way, successively more hat orientations can be found out - but this is very tedious, and can take a long time to find out for a single hat in the tiling, in some cases it is even necessary to first compute a larger part of the tiling, until the orientation can be determined.
\begin{figure}[bt]
	\begin{center}
		\includegraphics[width=\wwww]{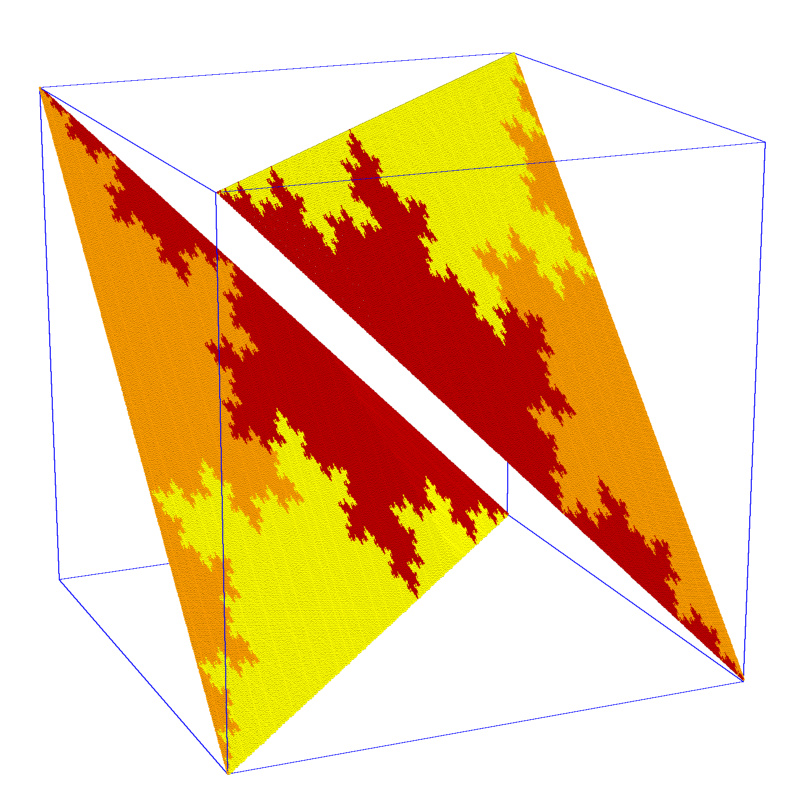}
		\hspace{0.01\textwidth}
		\includegraphics[width=\wwww]{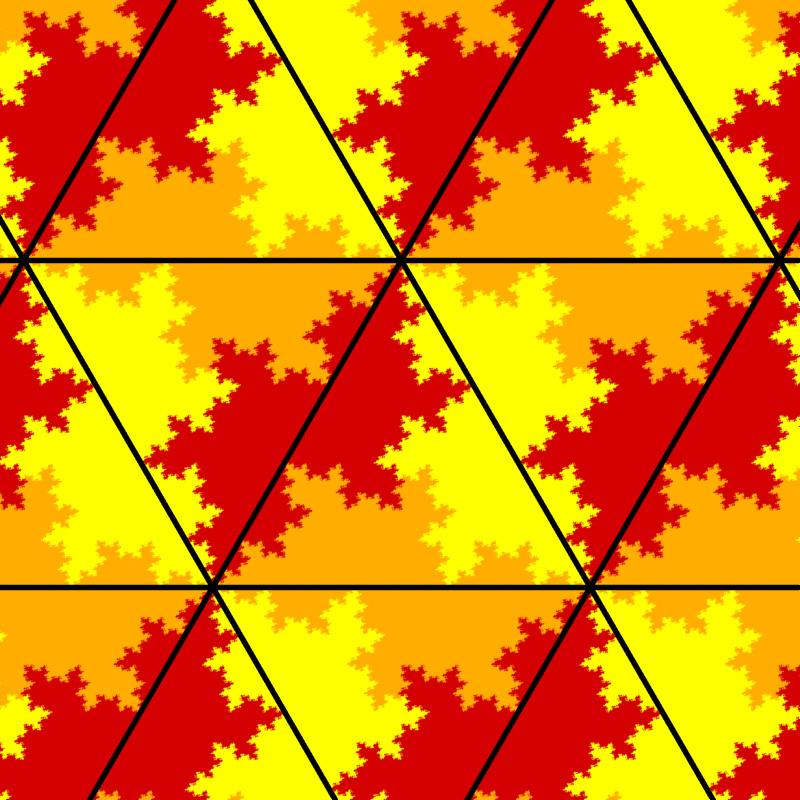}
		\hspace{0.01\textwidth}
		\includegraphics[width=\wwww]{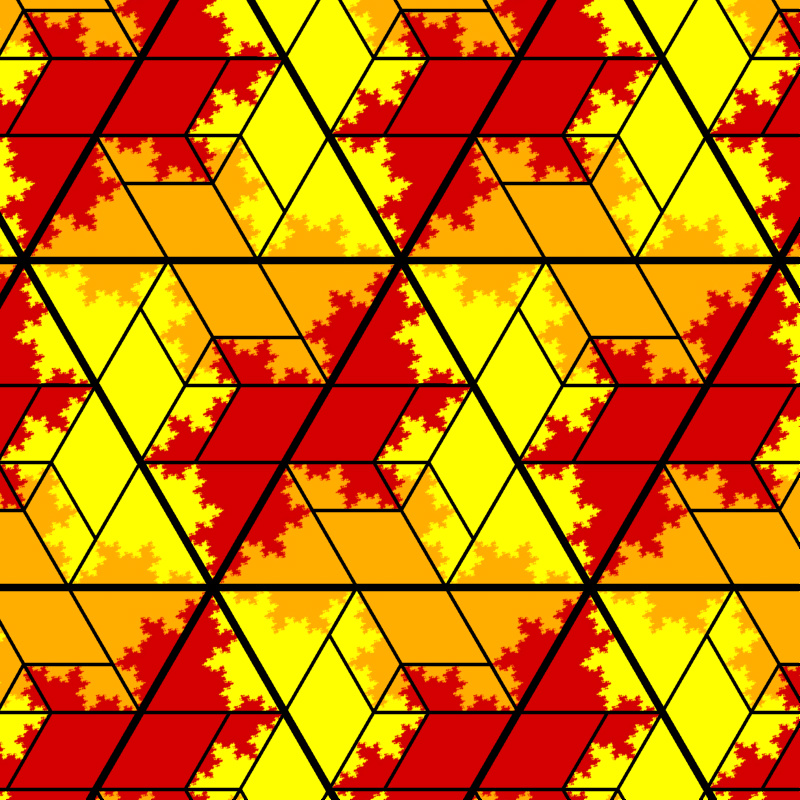}
		\hspace{0.01\textwidth}
		\includegraphics[width=\wwww]{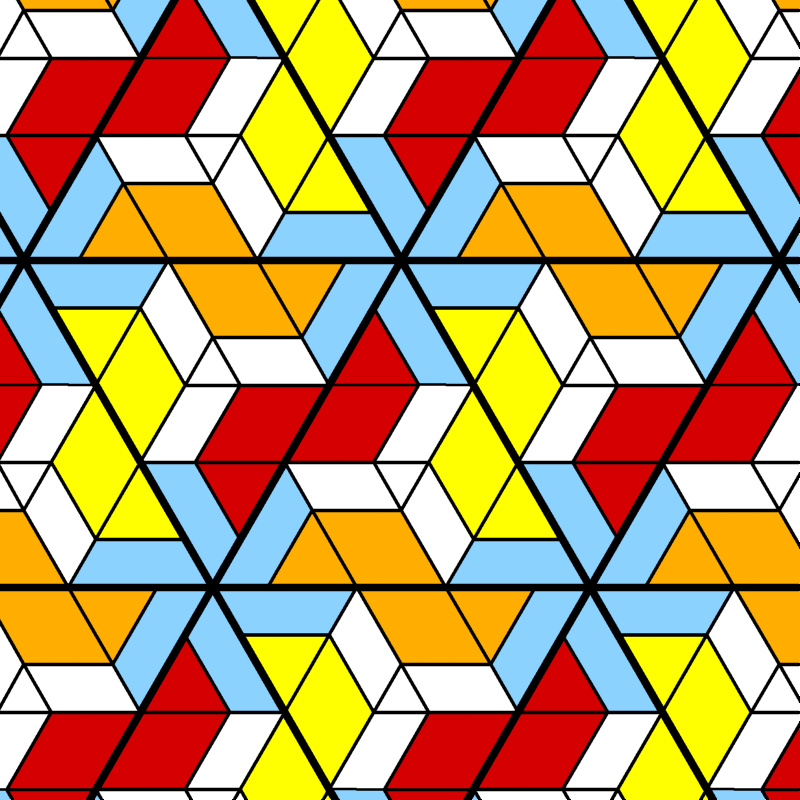}
		\caption{Left to right: Vertices $v_b$ coloured by hat orientation, showing a fractal structure -- Periodic orientation pattern in $P$ -- Orientation pattern in relation to the hat type pattern -- only vertices representing grey hats and first known light blue hat coloured by orientation.
		}	
		\label{orientationPattern}
	\end{center}
\end{figure}

Whereas the orientations for all grey hats are known in the first two steps, the orientations of light blue, dark blue and white hats are discovered successively along the tree-like  structure not covered by grey hats. At the isolated white hats -- that is where this structure branches -- the discovery of hat orientations has to be done from two directions until it can proceed into the third direction.

When we look at the orientations in the $v_b(i)$ representation of the triangle mesh, we observe a fractal-like structure, see Figure~\ref{orientationPattern}.

After colouring the pattern $P$ by hat orientation for the grey hat areas and the areas for the first light blue hats with known orientation, all areas for white hats and a part of the areas for light blue hats are still unknown.
To understand the fractal structure of the orientation pattern, we look at three structures in the hat tiling: a triple of light blue hats (neighbouring the same dark blue flipped hat), a region close to a pair of white hats, and the neighbouring six light blue hats to an isolated white hat.
By these three structures, we will be able to see local symmetries in the orientation pattern, which cover the missing area and help to compute orientations.

\begin{figure}[bth]
	\begin{center}
		\includegraphics[width=\www]{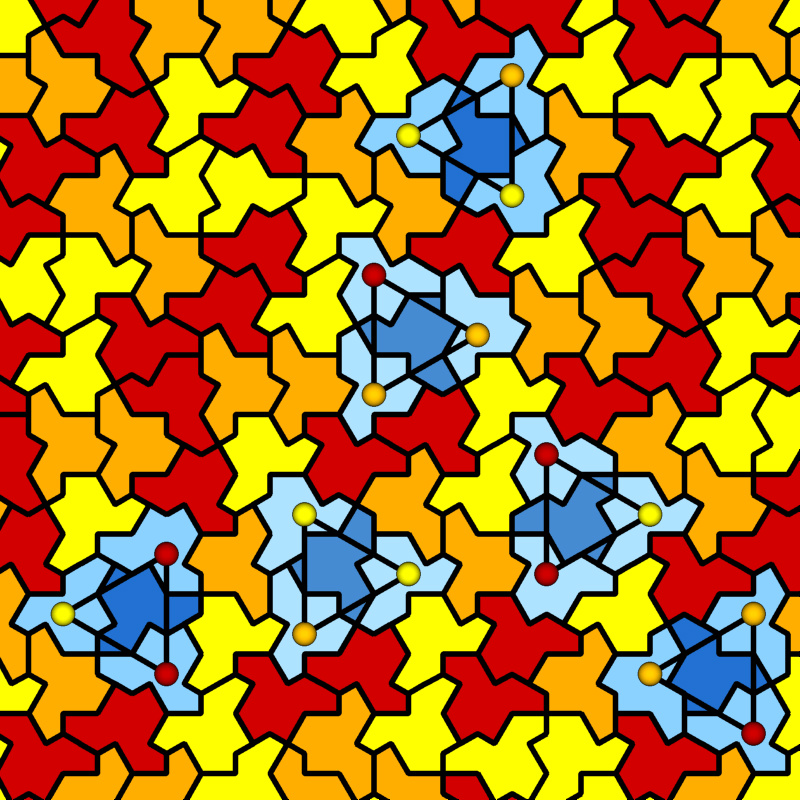}
		\hspace{0.01\textwidth}
		\includegraphics[width=\www]{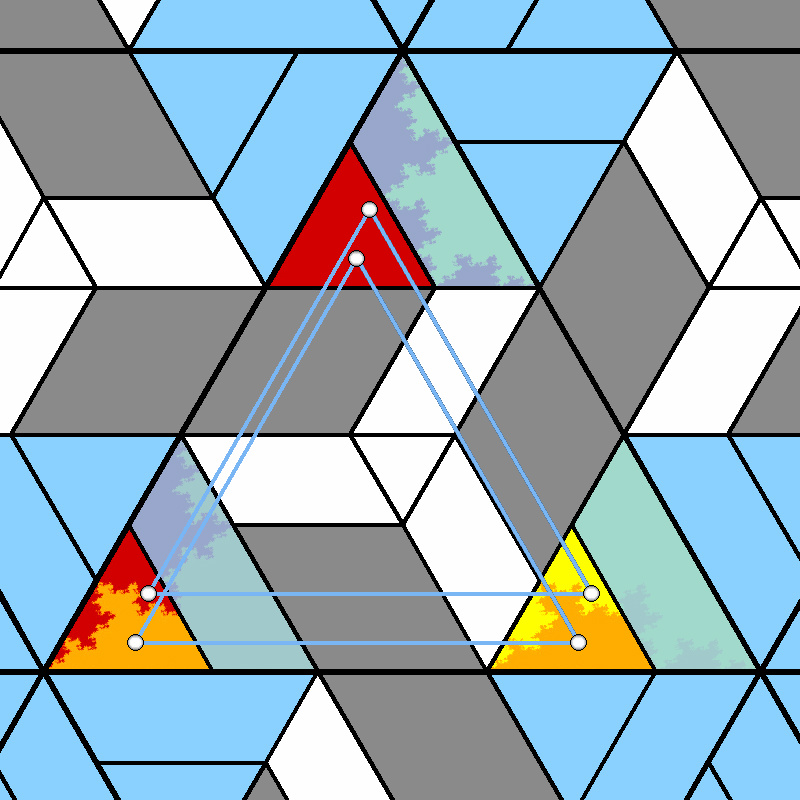}
		\hspace{0.01\textwidth}
		\includegraphics[width=\www]{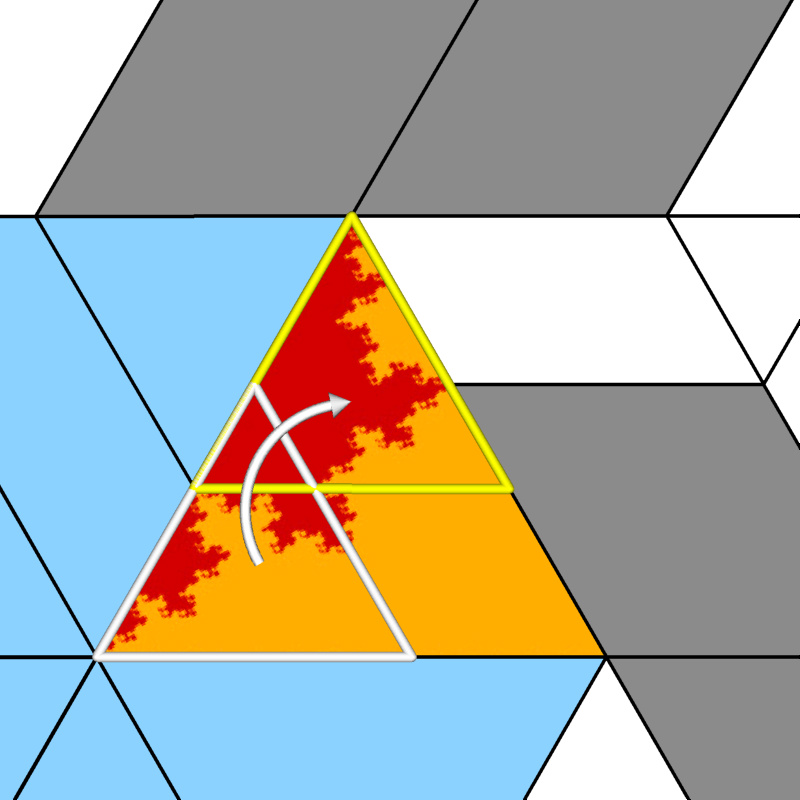}
		\caption{Left to right: Possible orientations of three light blue hats next to a flipped dark blue hat -- Equal pattern in small triangles -- Symmetry exchanging colours in one light blue area in $P$.
		}	
		\label{threeBlue}
	\end{center}
\end{figure}
Three light blue have six possible orientations as a group, three for centre triangles pointing up, and three for centre triangles pointing down.
If one of the three vertices of $T$ representing the three light blue hats is in the already known triangle area, there are two equally sized triangles in the other two light blue areas, and the pattern in one is a translation in the other triangle, with a cyclic permutation of orientation colours. Since the whole orientation pattern has a rotational symmetry by $\frac{2\pi}{3}$ (with cyclic permutation of orientations), in each of the three light blue triangles in $P$ there is one already completely known triangle with constant colour, and two triangles of equal size with yet unknown pattern, but there is a rotational symmetry from one to the other, exchanging the two orienations in this small triangle, see Figure~\ref{threeBlue}.

\begin{figure}[bt]
	\begin{center}
		\includegraphics[width=\www]{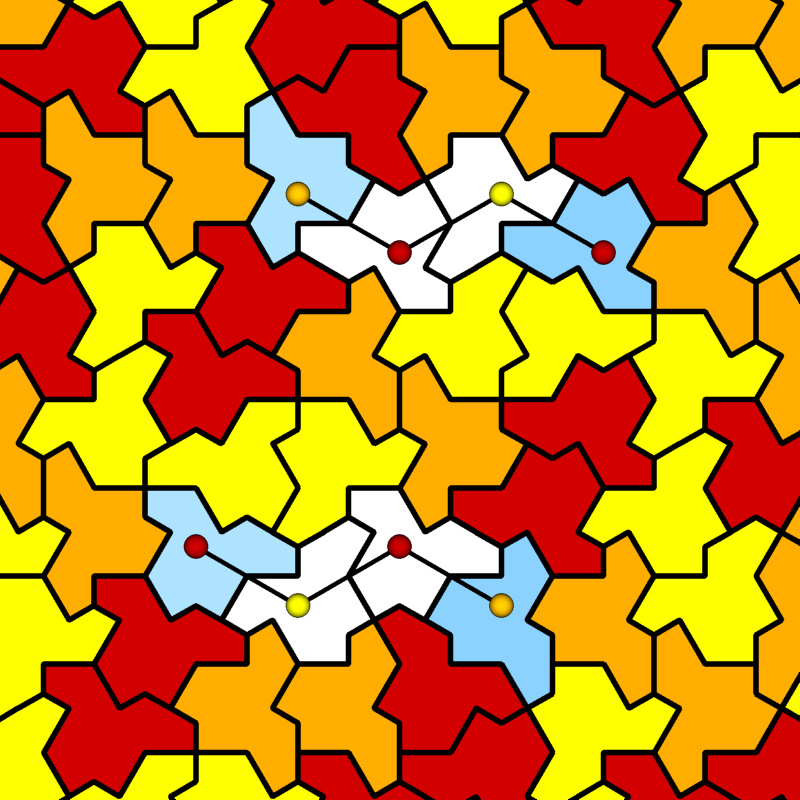}
		\hspace{0.01\textwidth}
		\includegraphics[width=\www]{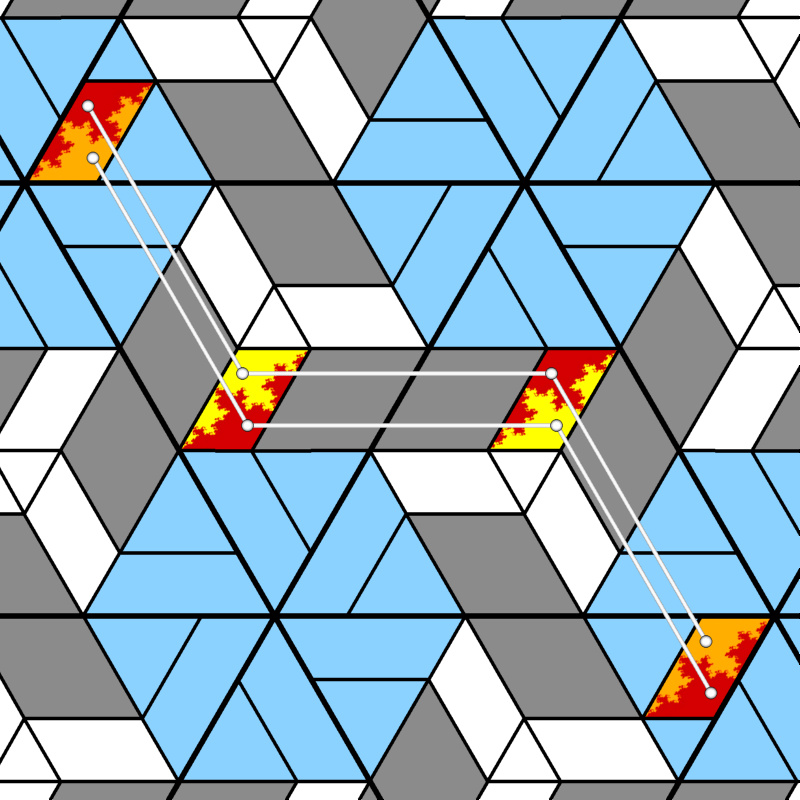}
		\hspace{0.01\textwidth}
		\includegraphics[width=\www]{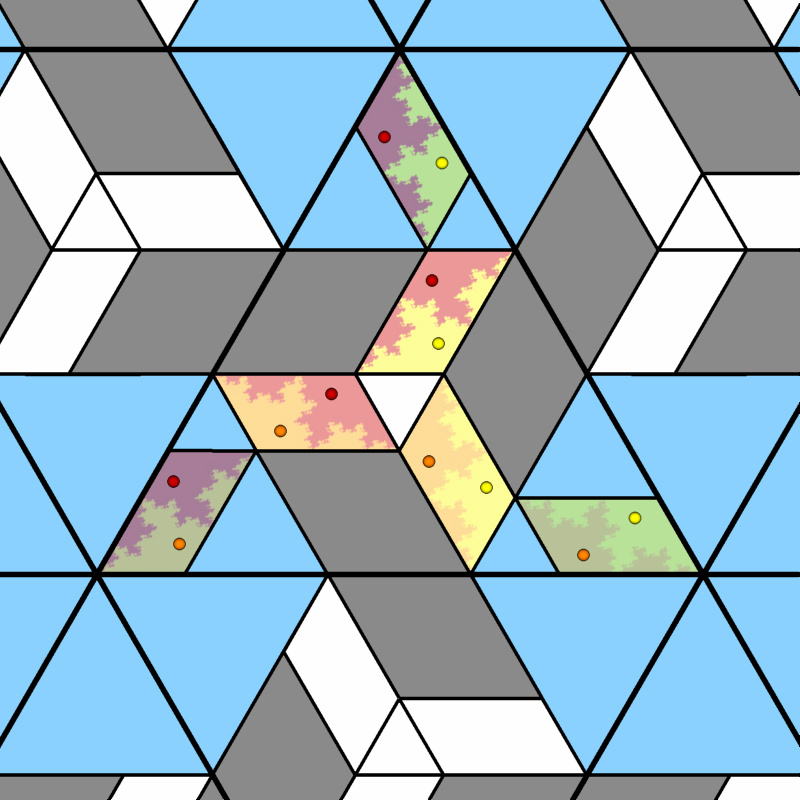}
		\caption{Left to right: Possible orientations for a pair of neighbouring white hats and two light blue hats - Equal pattern in parallelogram areas -- Symmetries in parallelogram areas in $P$.
		}	
		\label{whitePair}
	\end{center}
\end{figure}

By a pair of neighbouring white hats, two white parallelogram areas are related to each other. For a white pair with given positions, there are only two possible orientations, and these also determine the orientation of two light blue hats, see Figure~\ref{whitePair}.
Each of the shown parallelogram areas has rotational symmetry by $\pi$ to the centre of the parallelogram, exchanging the two occurring orientations (if it would not be symmetric, rotating the whole tiling by $\pi$ would require different orientations).
Once the orientation of one of the four hats is known, the other three orientations are fix as well.
So the parallelogram in the light blue area in $P$ has the same pattern as the white parallelogram in $P$ (again suitably permuting orientation colours). In one Triangle of $P$, there are six of these related parallelogram areas, each of them with an internal rotational symmetry by $\pi$.

\begin{figure}[bt]
	\begin{center}
		\includegraphics[width=\www]{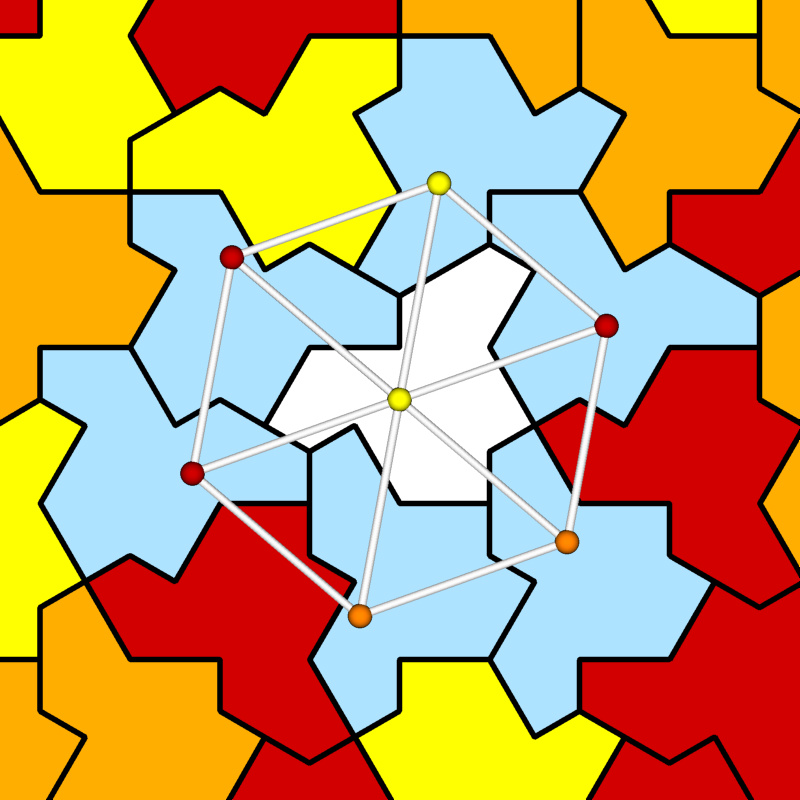}
		\hspace{0.01\textwidth}
		\includegraphics[width=\www]{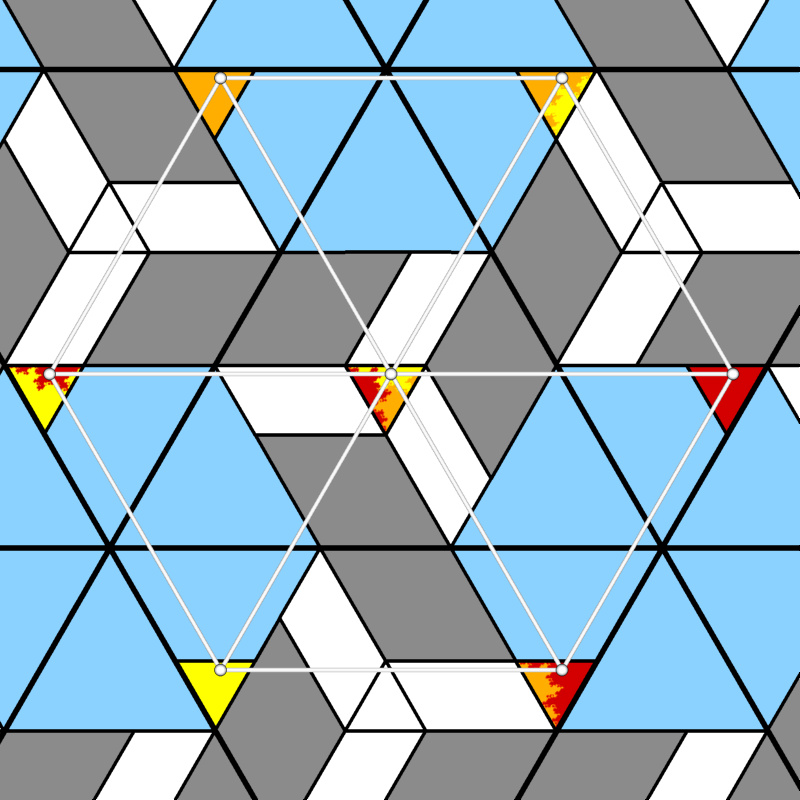}
		\hspace{0.01\textwidth}
		\includegraphics[width=\www]{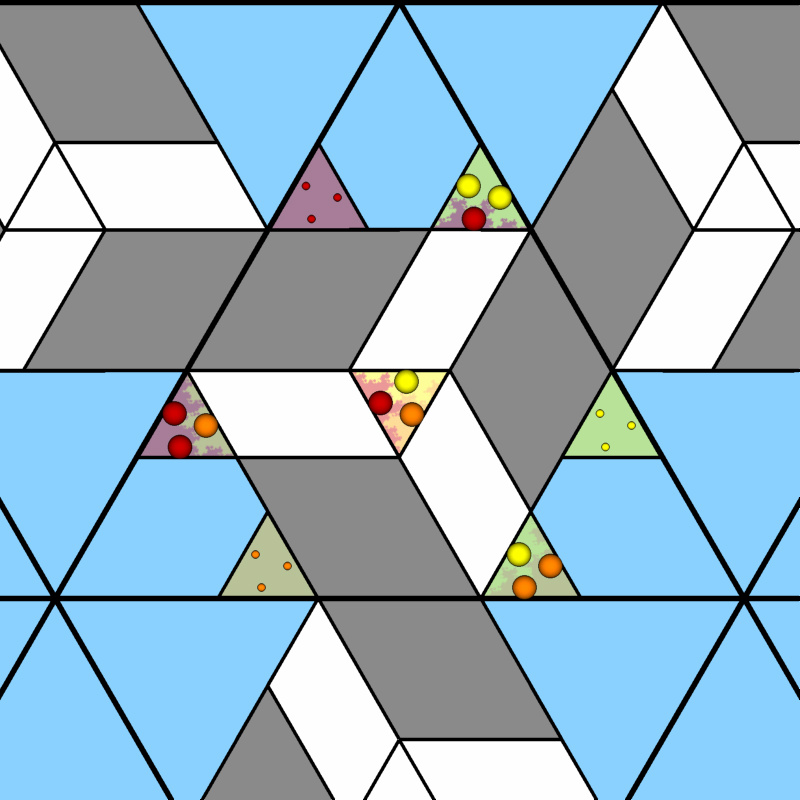}
		\caption{Left to right: Possible orientation for a the neighbourhood of an isolated white hat - Equal pattern in small triangle areas -- Symmetries of small triangles in the orientation pattern in one triangle of $P$.
		}	
		\label{isolatedwhite}
	\end{center}
\end{figure}

The two discussed partial symmetries in the orientation pattern do not yet include the centre white triangle in $P$; this centre triangle represents isolated white hats.
Isolated white hats have a neighbourhood of six light blue hats, and there are six possible combinations of orientations for the whole group - three for the centre triangle of the isolated white hat pointing up, three for centre triangle pointing down. One of the possible neighbourhoods is shown in Figure~\ref{isolatedwhite}.
The orientations of three of the neighbouring light blue hats are independent of the relative position of the vertex representing the isolated white hat, the other three draw a similar pattern as in the centre white triangle into three small triangles of the same size in $P$ into light blue areas, but using only two of the three possible orientations:
if the vertex in the centre triangle and the vertices representing neighbour light blue hats in the unknown area are drawn into the same triangle of the pattern $P$, and eight more vertices are drawn into the same triangle of $P$ by $\frac{2\pi}{3}$ rotational symmetry, cyclically permuting the orientation colours, we get a set of twelve related vertices.
In the centre triangle, all three colours are used, the other three triangles have two vertices of one colour and one vertex of a different colour.
Once one colour in the centre triangle is known, all other vertices can be coloured as well.
In reverse, if for two $\frac{2\pi}{3}$ rotationally symmetric areas in one of the small triangles not in the centre the colours are known, the third symmetric area in the same triangle and all nine symmetry images in the other three small triangles are known as well.

By the three discussed symmetry operations, the whole unknown area from the beginning can be discovered. Essentially these symmetry operations do the same as successively completing the six-kites parts of hats to complete hats once there is only one possibility left.
The advantage of the evaluation of the fractal orientation pattern is, that it can be evaluated for each hat individually. It is not necessary to grow the information along the white and blue tree structure, which in some cases requires the computation of a bigger part of the hat tiling.
In addition, in the iterative evaluation of the fractal pattern the unknown area shrinks exponentially with every iteration step, so this evaluation is faster than linearly following the branches in the white and blue tree structure in the hat tiling.

Because of the self--similarity of the fractal structure, we can evaluate the orientation in an unknown area iteratively as follows:
\begin{itemize}
\item If the position is in the centre white triangle of $P$, scale the centre triangle to the complete $P$ triangle and rotate by $\frac{\pi}{3}$.
\item If the position is in a white parallelogram area, we scale the parallelogram by $\Phi^2$.
\item If the position is in a light blue area, we scale by $\Phi^2$.
\end{itemize}
\begin{figure}[bt]
	\begin{center}
		\includegraphics[width=\www]{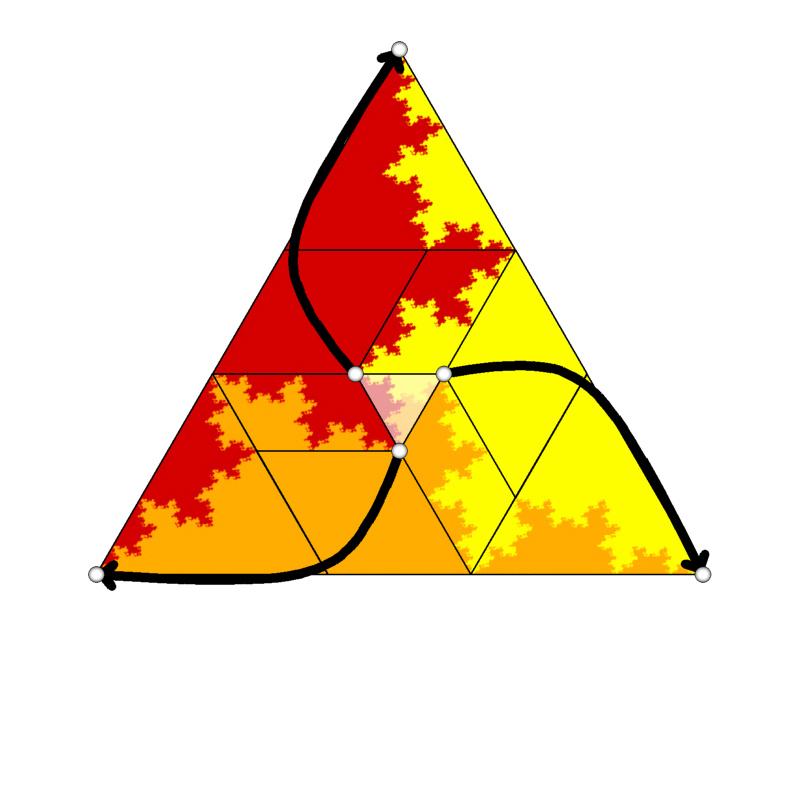}
		\hspace{0.01\textwidth}
		\includegraphics[width=\www]{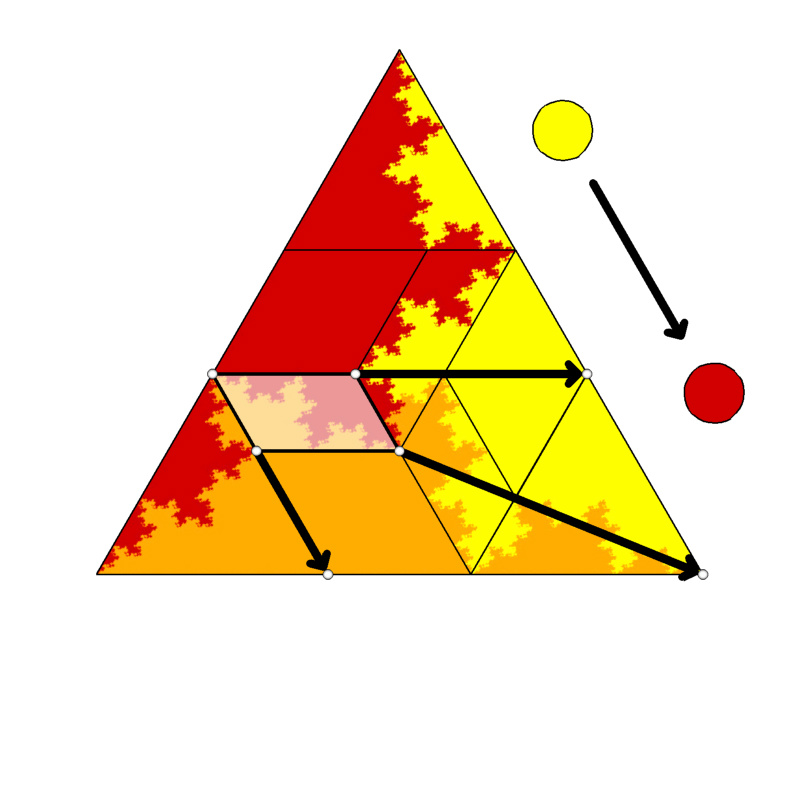}
		\hspace{0.01\textwidth}
		\includegraphics[width=\www]{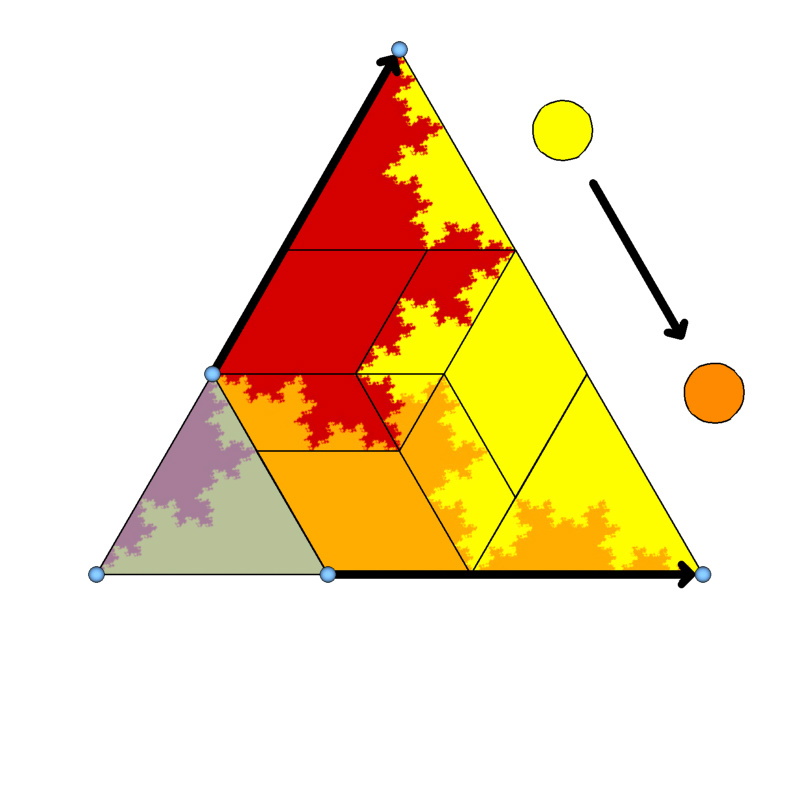}
		\caption{Left to right: Map the centre triangle to the whole $P$ triangle -- Scale unknown white parallelogram area and change colour -- Scale light blue triangle area and change colour.
		}	
		\label{evaluateFractal}
	\end{center}
\end{figure}
If at some time in this iterative process the point was in a white parallelogram or light blue triangle, we replace one of the colours of the orientation pattern in $P$ and use only the two appropriate colours used in the smaller area, see Figure~\ref{evaluateFractal}.
Since with the number of iterations the unknown area shrinks exponentially, the orientation for a position is in most cases discovered by a very small number of iterations.

\section{Construction summary}
So, to construct a hat tiling, we perform the following steps:
\begin{enumerate}
\item We chose two of the $d_k$ values, the third is then known by $\sum d_k = 0$.
\item We compute integer indices $a_i$, $b_i$ and $c_i$ for all vertices $v_i$ in $T$.
\item From $a_k(i)$ and $b_k(i)$ we evaluate the $v_b(i)$ positions.
\item From $v_b(i)$ and $P$ we evaluate the hat type.
\item From $v_b(i)$ and the fractal orientation pattern we evaluate the hat orientation.
\item From $c_k(i)$, the position of the hat centre triangle is known - we draw a hat tile with the correct orientation around the centre triangle.
\end{enumerate}

\section{Flipping hats}
\begin{figure}[bt]
	\begin{center}
		\includegraphics[width=\wwww]{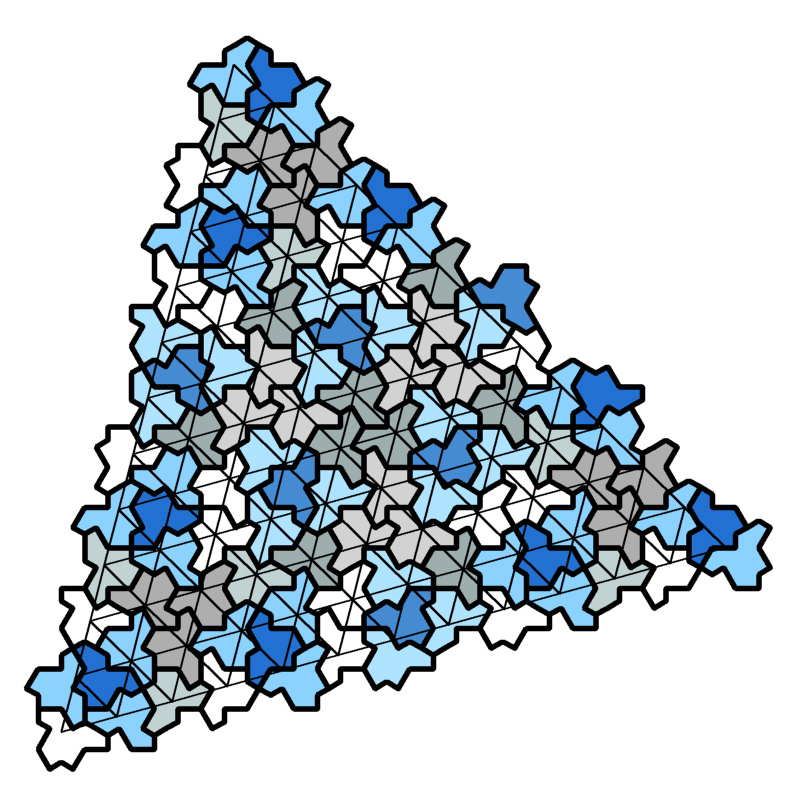}
		\hspace{0.01\textwidth}
		\includegraphics[width=\wwww]{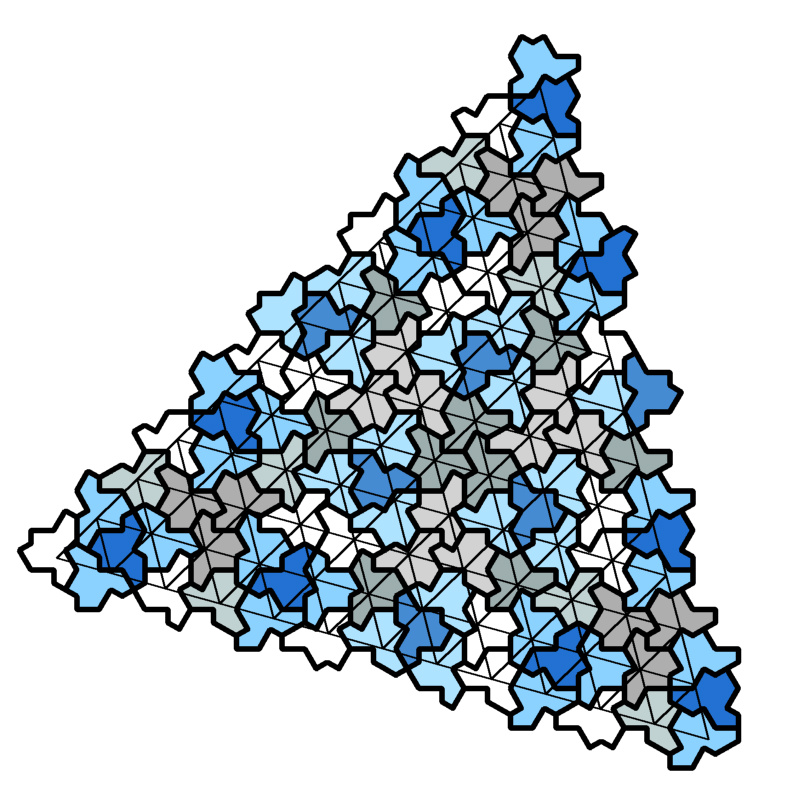}
		\hspace{0.01\textwidth}
		\includegraphics[width=\wwww]{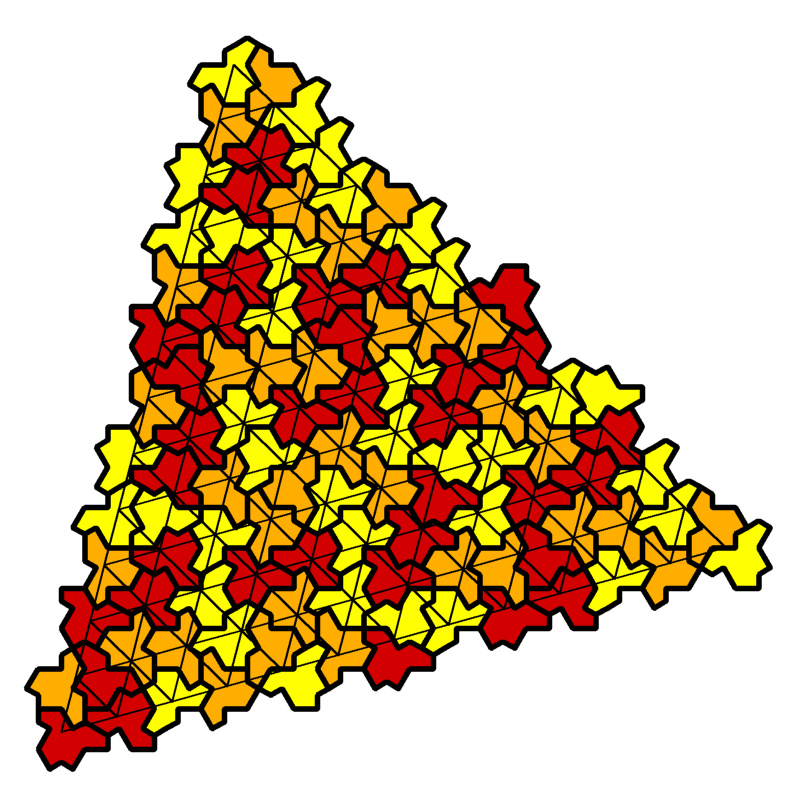}
		\hspace{0.01\textwidth}
		\includegraphics[width=\wwww]{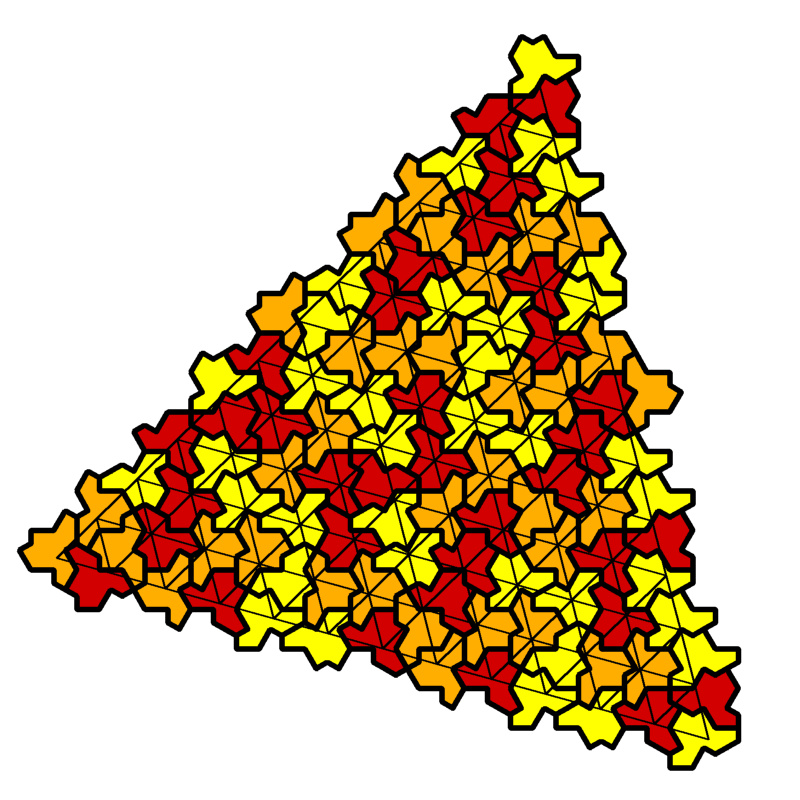}
		\caption{Same section of a hat tiling with exchanged roles of flipped and unflipped hats, coloured by hat types and orientations.
		}	
		\label{flipHats}
	\end{center}
\end{figure}
When we exchange the roles of flipped and unflipped hats, i.e. both the flipped and unflipped hat are replaced by their mirror images, a pattern can be computed from the same data; the computation of $a_k(i)$ and $b_k(i)$ does not change.
For $c_k[i]$, we use Equation~\ref{c2} instead of equation~\ref{c}. In the computation of hat types and orientations, we simply have to use the mirror reflection of $P$.
With the exchanged roles, the light blue and isolated white hats keep their type, but some grey hats become white and some white hats become grey.
The fractal structure of the orientation pattern makes the changes of orientations difficult to see.
We can see, that dual trianglulation $T$ is rotated in different directions for exchanged roles of flipped and unflipped hats, see Figure~\ref{flipHats}.
\begin{figure}[bt]
	\begin{center}
		\includegraphics[width=\wwww]{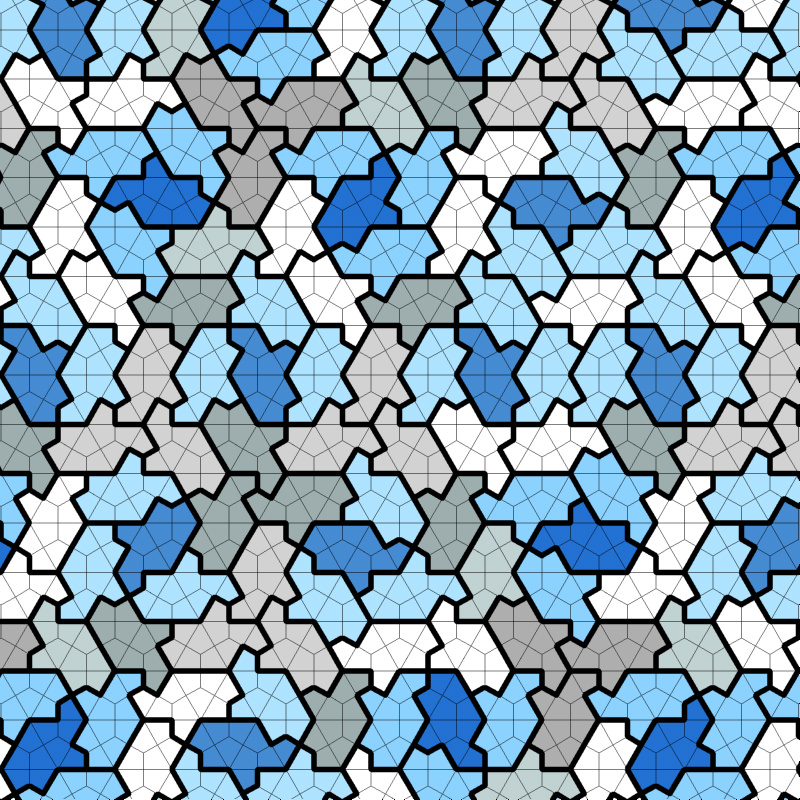}
		\hspace{0.01\textwidth}
		\includegraphics[width=\wwww]{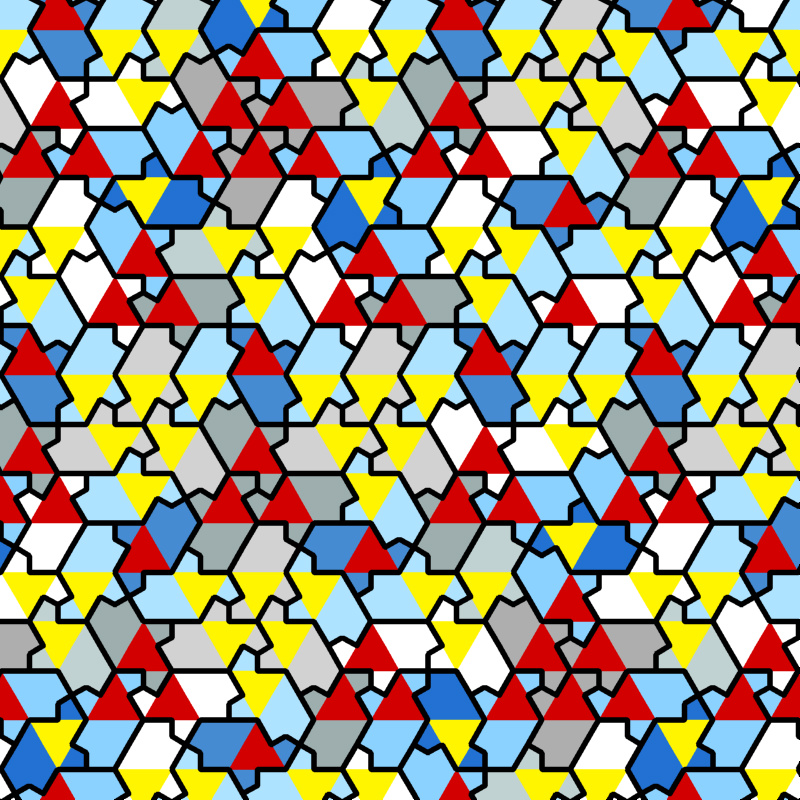}
		\hspace{0.01\textwidth}
		\includegraphics[width=\wwww]{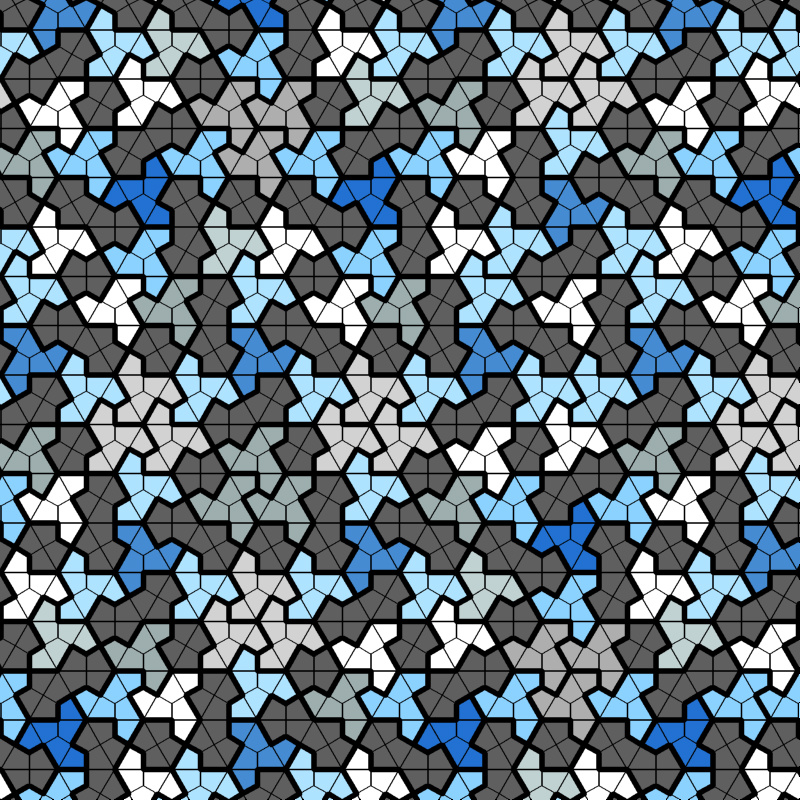}
		\hspace{0.01\textwidth}
		\includegraphics[width=\wwww]{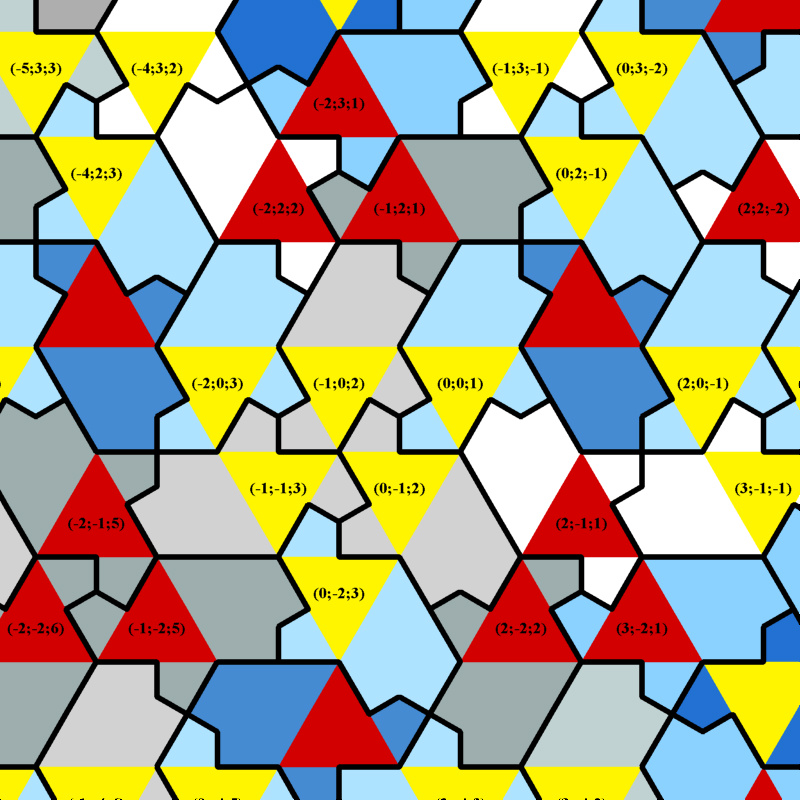}
		\caption{$10$--kites Tile; left to right: Part of a tiling -- the centre triangles, red pointing up, yellow pointing down -- rotational symmetric part consisting of six kites in each tile -- $c_k$ labels indicating positions centre triangles.
		}	
		\label{center10}
	\end{center}
\end{figure}

For the $10$--kites monotile, the dual triangulation connecting the centre triangles is not rotated in relation to $U$. The two versions with exchanged roles of flipped and unflipped $10$--kite tiles have the same realisation of $T$ with vertices in the unflipped hats.
For the $10$--kites tile, the positions for centre triangles can be computed by
\begin{align}
c(i) = \begin{pmatrix} c_0(i) \\ c_1(i) \\ c_2(i) \end{pmatrix} = \begin{pmatrix}a_0(i)+2 b_0(i) \\ a_1(i)+2 b_1(i) \\ a_2(i)+2 b_2(i)\end{pmatrix},\label{c3}
\end{align}
in this case the formula does not change for exchanged roles of flipped and unflipped tiles. The $10$--kites tile contains two complete triangles from $U$, here the triangle with three rotational symmetric kites araound belonging to the tile is considered the centre triangle, see Figure~\ref{center10}.

\begin{figure}[b]
	\begin{center}
		\includegraphics[width=\wwww]{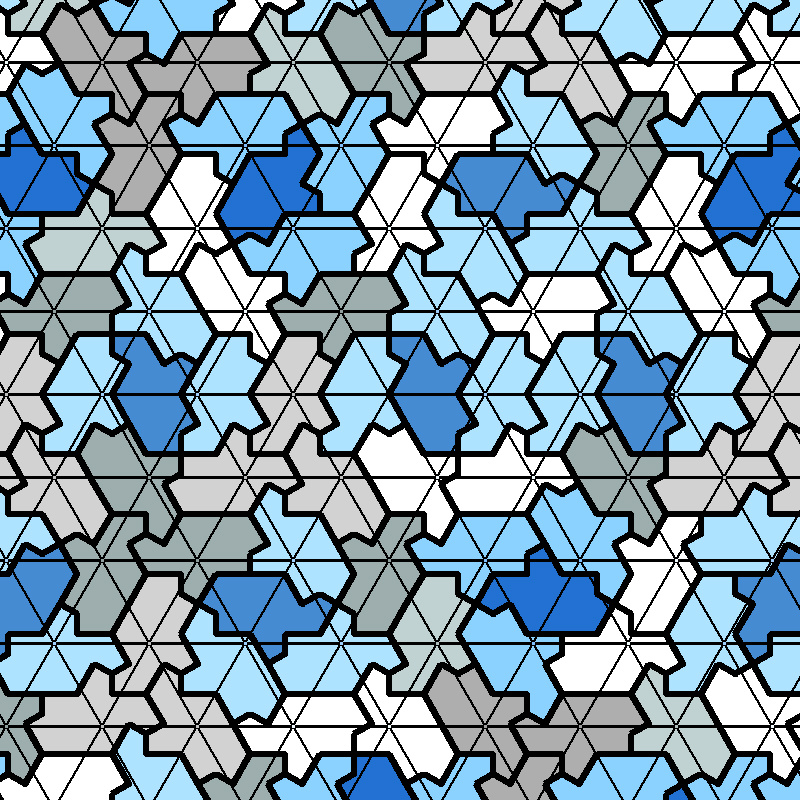}
		\hspace{0.01\textwidth}
		\includegraphics[width=\wwww]{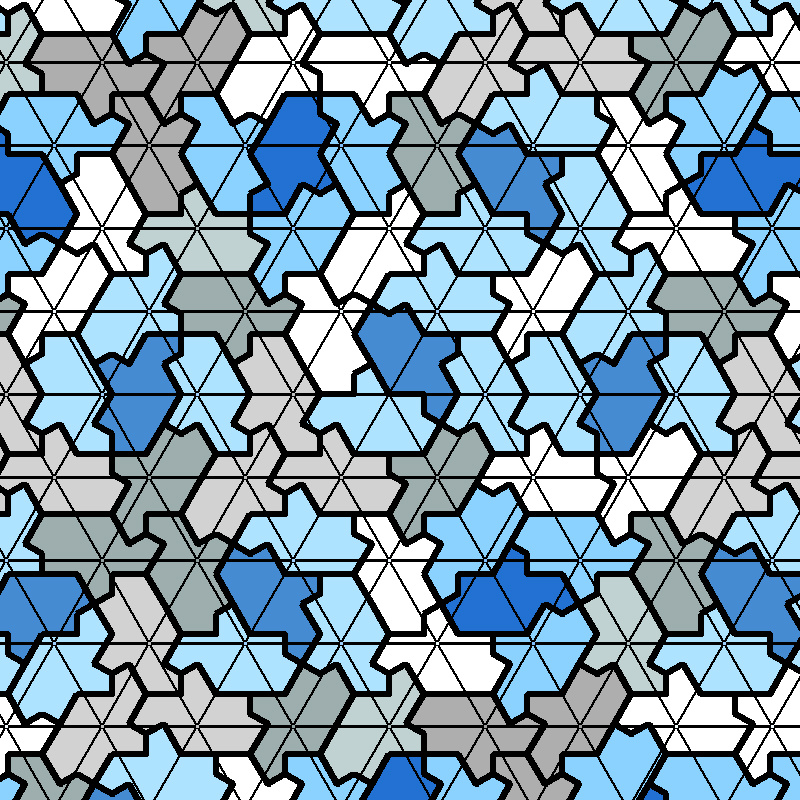}
		\hspace{0.01\textwidth}
		\includegraphics[width=\wwww]{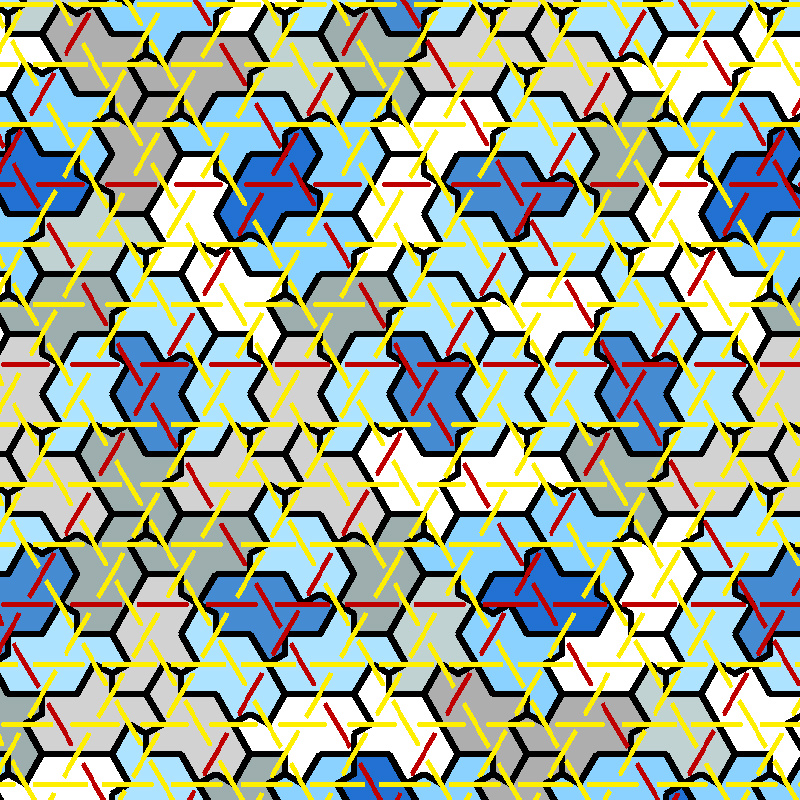}
		\hspace{0.01\textwidth}
		\includegraphics[width=\wwww]{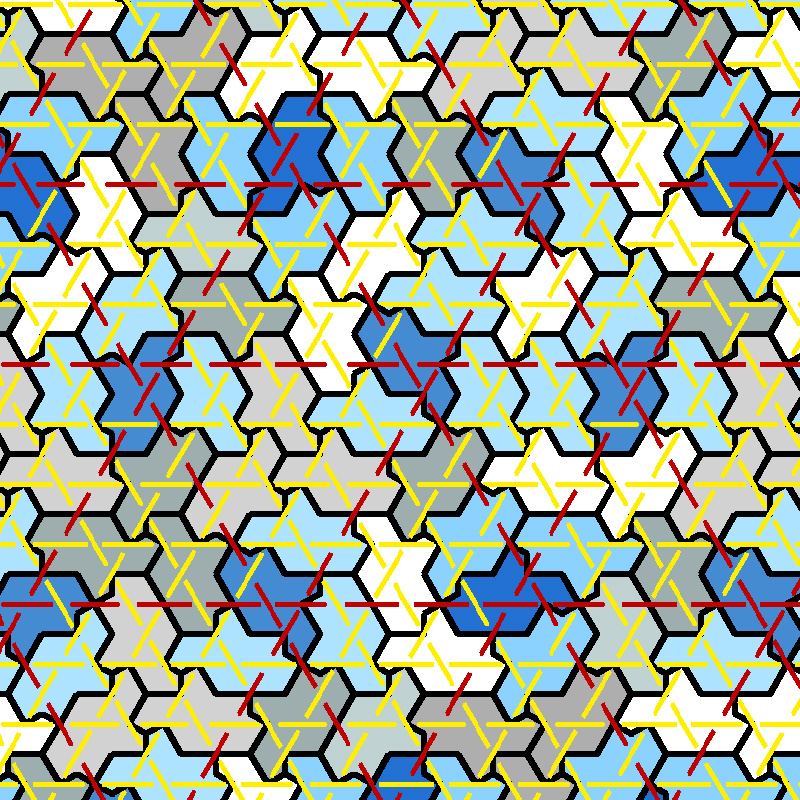}
		\caption{Same section of a $10$--kite tiling with exchanged roles of flipped and unflipped hats, coloured centre lines are the same for both versions.
		}	
		\label{tenKites}
	\end{center}
\end{figure}

When we draw centre lines between neighbouring lines in $U$, some of these pass through centre triangles of unflipped tiles only, the other ones pass through centre triangles of flipped tiles only.
Colouring these two types of centre lines differently, each unflipped $10$--kite tile shows the same line decoration, and the flipped ones show the same line decoration with exchanged colours, see Figure~\ref{tenKites}.
So the $10$--kites tile is in some sense less dependent on the choice which tile is the flipped and which is the unflipped type.


\bibliographystyle{plain}
\bibliography{monotile} 
\end{document}